\documentclass[11pt]{amsart}
\usepackage{amssymb,amsthm,amsmath,epsfig}
\usepackage{calc}
\usepackage{times}

\newcommand{\rar}{\rightarrow}
\newcommand{\lar}{\longrightarrow}
\newcommand{\surjects}{\twoheadrightarrow}

\newcommand{\C}{{\mathbb C}}

\newcommand{\fm}{{\mathfrak m}}
\newtheorem{defn0}{Definition}[section]
\newtheorem{prop0}[defn0]{Proposition}
\newtheorem{conj0}[defn0]{Conjecture}
\newtheorem{thm0}[defn0]{Theorem}
\newtheorem{lem0}[defn0]{Lemma}
\newtheorem{corollary0}[defn0]{Corollary}
\newtheorem{example0}[defn0]{Example}
\newtheorem{quest0}[defn0]{Question}
\newtheorem{rmk0}[defn0]{Remark}

\newenvironment{defn}{\begin{defn0}}{\end{defn0}}
\newenvironment{prop}{\begin{prop0}}{\end{prop0}}

\newenvironment{thm}{\begin{thm0}}{\end{thm0}}
\newenvironment{lem}{\begin{lem0}}{\end{lem0}}
\newenvironment{cor}{\begin{corollary0}}{\end{corollary0}}
\newenvironment{exm}{\begin{example0}\rm}{\end{example0}}
\newenvironment{quest}{\begin{quest0}\rm}{\end{quest0}}
\newenvironment{rmk}{\begin{rmk0}\rm}{\end{rmk0}}

\newcommand{\demo}{{\bf Proof. }}

\begin{document}

\title {Homology of homogeneous divisors}

\author{Aron Simis}
\thanks{The first author is partially supported by a CNPq Grant.}
\address{Universidade Federal de Pernambuco\\ Departamento de Matemática\\ Av. Prof. Luiz Freire\\ s/n
Cidade Universitária, CEP 50740-560\\ Recife - Pernambuco - Brasil\\}
\email{aron@dmat.ufpe.br}

\author{\c Stefan O. Toh\v aneanu}
\address{Department of Mathematics\\ The University of Western Ontario\\ London, Ontario N6A 5B7\\ Canada\\}
\email{stohanea@uwo.ca}


\begin{abstract}
\noindent One deals with arbitrary reduced free divisors in a polynomial ring over a field of characteristic
zero, by stressing the ideal theoretic and homological behavior of the corresponding singular locus.
A particular emphasis is given to both weighted homogeneous and homogeneous polynomials, allowing to introduce
new families of free divisors which do not come from hyperplane arrangements nor as explicit discriminants from
singularity theory.
\end{abstract}
\maketitle


\section*{Introduction}

Let $Y$ be a reduced divisor on a smooth algebraic variety $X$ over the complex field. According to \cite{sa},
$Y$ is called a {\em free divisor}
on $X$ if the $\mathcal O_X-$module $${\rm Der}_X(-\log Y):=\{\theta\in {\rm Der}(X)\,|\,\theta(\mathcal O_X(-Y))
\subseteq \mathcal O_X(-Y)\},$$
is free, where $\mathcal O_X$ is the sheaf of regular functions on $X$.

A special case is that of a non-smooth reduced divisor $Y=V(F)\subset\mathbb P^2$ on the projective plane, where
$F$ is a squarefree homogeneous
form of degree $d$ in $R:=\mathbb K[x,y,z]$.
$F$ is called {\em the defining polynomial of $Y$} and the ideal $J_F$ of $R$ generated by the partial derivatives
$F_x, F_y, F_z$ of $F$
is the {\em gradient ideal} of $F$.
Since $F$ is squarefree and non-smooth, $J_F$ has codimension $2$.
Let ${\rm Der}(R)$ be the $R-$module of derivations on $R$, and let $$D_F:=\{\theta\in {\rm Der}(R):\theta(F)\in \langle F\rangle\}$$ be
the $R-$module of logarithmic derivations on $F$. $Y$ is called a {\em free divisor} if and only if $D_F$ is a free $R-$module.
Since $dF=xF_x+yF_y+zF_z$, one has that $\theta_E=x\frac{\partial}{\partial x}+y\frac{\partial}{\partial y}+z\frac{\partial}{\partial z}$
is in $D_F$. In fact, $$D_F=\theta_ER\oplus D^0_F,$$ where $D^0_F$ is an $R-$submodule of $D_F$ whose elements are
in one-to-one correspondence with the syzygies on $J_F$.
We have that $Y$ is free if and only if $J_F$ is a perfect ideal (i.e., the ring $R/J_F$ is Cohen-Macaulay).

The notion of a free divisor was originally associated with hyperplane arrangement theory.
When $Y$ is a hyperplane arrangement in a vector space of finite dimension, free divisors received a great deal of interest:
\cite{ot,sa,sc,te1,te2,y1} to cite just a few. When $X$ is not necessarily the projective or the affine space,
\cite{dssww} gives some insights on the freeness of divisors on $X$ in this general setup. In the case of divisors
on the projective space, or affine space, some important results were obtained, especially when $X=\mathbb P^2$ in
\cite{AA,st,s1,s3,t1}. The advantage of working with divisors on $\mathbb P^n$ is that Saito's Criterion (\cite{sa})
translates into the following: the divisor $Y$ is free if and only if the gradient ideal of $Y$ is a perfect ideal of codimension 2.
On $\mathbb P^2$, because this ideal is generated by three homogeneous polynomials of the same degree, ideals of this type
have been studied extensively (besides the already mentioned references, see \cite{cs,HS,SV1}). The disadvantage is that
for these divisors one does not have a lattice of intersection, and hence the combinatorics plays a very
minimal role in the game.

The tone of free divisor theory as embedded in hyperplane arrangements is that the divisors come up naturally with
many irreducible components.
Discriminant theory gives another picture of the theory away from arrangements, but a systematization thereof is yet to come
up in the literature.
Even when keeping offshore from hyperplane arrangements, the results have remained close to generalizations thereof,
such as arrangements of lines and smooth conics in $\mathbb P^2$ (as in \cite{st}).
A question arises as to whether there are families of homogeneous free divisors of arbitrary degrees which arise out
of simple constructs in algebraic geometry.
This question makes sense even for non-irreducible reduced homogeneous polynomials, but is reaches its pick interest
in the irreducible case.

By a suitable {\em addition} procedure one can turn some highly nonfree irreducible homogeneous divisors into free divisors
(see Proposition~\ref{addition} and Proposition~\ref{addition2}).
As an illustration, in $3$ variables and arbitrary degree, this procedure yields a free divisor  $f=gh$ with just two irreducible
factors, where $g=y^rz^{d-r}-x^d$ is a cuspidal
type singularity and $h=y$ is its multiple tangent line at $(0:0:1)$.

 At the other end, there are homogeneous {\em irreducible} free divisors of linear type of degree $4$ in $4$ variables.
One such example is the defining equation of the tangential surface (or the dual surface) to the rational normal
curve in $\mathbb{P}^3$ (\cite[Remark 4.4]{s1}) -- it is highly suggestive that this is the case of the dual hypersurface
$V(F)\subset (\mathbb{P}^n)^*$ of the rational normal curve $C_n\subset \mathbb{P}^n$, for any $n\geq 3$, an
expectation that has been verified computationally
for initial values of $n$. Alas, the degree of $V(f)$ is $2n-2$, so we get no examples in odd degrees.
Thus, even in higher dimension, homogeneous irreducible free divisors are not so abundant, as it seems.
But it is in dimension $2$ that things get really tight, thus leaving the hope to be able
to classify them. In most of the papers cited here, to be able to say something about
the freeness of these ``nonlinear'' divisors, an extra condition about the singularities must be added: they must be
locally everywhere weighted homogeneous.
In terms of the gradient ideal, this ideal must be locally (at the minimal primes) a complete intersection (see \cite{r} or,
for some more detailed explanations, \cite{st}).

The overall goal of this work is an insight into the freeness of divisors on $\mathbb P^n$,
with particular emphasis in the case $n=2$.
Some of the main results are stated in Theorem~\ref{Cramer_vs_sat}, Theorem~\ref{three_vars} and Theorem~\ref{addition2}.

Let us now briefly describe the contents of each section.

The first section is a collection of known and new results about codimension $2$
ideals $I$ in the polynomial ring $R:=\mathbb K[x,y,z]$ which are generated by three forms (almost complete intersection) of same degree.
The goal is to have a priori properties enjoyed by the singular ideal (gradient ideal) of a homogeneous
polynomial and to establish possible templates for the latter to fit in.
These properties distant themselves from combinatorics in that they mainly appeal to the usual numerical
invariants coming from commutative algebra and homology.
In this line of quest, we discuss the syzygies and the regularity of an ideal of the above kind.
The behavior of the latter helps to distinguish between a perfect such ideal (the case of the gradient ideal
of a free divisor) and a non-perfect one.
Some of these results are in  line with with \cite{HS}.
We then consider the symmetric and Rees algebras of the ideal since these configure a way of understanding
when a free divisor is Koszul-free or of linear type, or even of syzygetic type.
There is a deep entanglement between these pre-homological properties and the geometry of the divisor which is yet to
be better understood.
The last part of the first section is mainly devoted to analyzing the local cohomology module $H_{\langle x,y,x \rangle}^0(R/I)=
I^{\rm sat}/I$ in the special case in which it coincides with the socle of $R/I$.
This led us naturally to introduce a class of ideals we call {\em Cramer ideals}.
Of course,  the whole section is moved by the spectrum of the following questions:
given a $3\times 2$ homogeneous  matrix $\phi$  over $R$ such that the ideal $I_2(\phi)$ of $2$-minors has codimension $2$,
when are these minors {\em integrable} (i.e., the partial derivatives of a form in $R$)?
For which non-perfect homogeneous ideals $J\subset R$ (not necessarily generated in same degree) is $J^{\rm sat}$ the gradient ideal
of a free divisor?

The second section is totally devoted to a discussion of free divisors.
We first get a glimpse  of irreducible homogenous free divisors according to their degree,
describing infinite families of such divisors of degree $\geq 5$.
Unfortunately, for any degree $\geq 6$ the gradient ideal of a member fails for the property of being
of linear type.
It is as yet a challenge to exhibit explicit families of irreducible homogeneous free divisors of linear type for all degree $\geq 6$.
In the sequel, we study weighted homogeneous polynomials and their homogenization as regards to  freeness.
We show that the gradient ideal of the homogenization always has a linear syzygy, a result that turns out to be very basic
in this part. Using it we are able to characterize the binary weighted homogenous divisors $f\in\mathbb K[x,y]$
whose homogenization $F\in \mathbb K[x,y,z]$ is a free divisor\,; in particular, no irreducible such $f$ has this property.
We are likewise able to write down the minimal free resolution of $J_F$ in any case and note that the corresponding
shifts depend only on $\deg(f)$ (not on the integer weights of $f$).
We also discuss the so-called {\em cone} cone$(f)\in \mathbb K[x_1,\ldots,x_n,t]$ of a weighted homogeneous
$f\in \mathbb K[x_1,\ldots,x_n]$ and give a sufficient condition for the cone of $f$ to be a free divisor in terms
of the homogenization of $f$.
The rest of the section looks at the details of some constructs, some new, some known.
Here we introduce some classes of homogenous  divisors which give rise to free homogenous divisors by means
of {\em adding} a free divisor.
The procedure yields an iterative method to create new free divisors of arbitrary degrees asymptotically
as the number of variables grows.
Finally, we review some classical free divisors in regards to the ideal theoretic and homological tool
introduced in earlier parts.

\section{Ideal theoretic and homological preliminaries}

Let $\mathbb K$ stand for an algebraically closed field of characteristic zero
(usually the complex numbers).
The restriction on the characteristic may be lifted for much of the preparatory material in this part.

\subsection{Syzygies and regularity.}
\label{ssec2.1}

Almost all the results in this subsection are valid for more general classes of rings and ideals.
Since our primary interests are codimension 2 ideals $I\subset R:=\mathbb K[x,y,z]$ minimally generated by three
homogeneous polynomials of the same degree, we will state the results accordingly to this setup.

We will need a preliminary result of general interest. It will suffice to consider a restricted case though is
admits more general formulations.
Its contents follow essentially from a close inspection of a Hilbert series (see \cite[Proof of Theorem 1.5, (13)]{HS}),
but we chose to give a completely elementary proof that also clarifies its use.

\begin{lem}\label{cod2_3gens}
Let $I\subset R=\mathbb K[x,y,z]$ be an ideal of codimension $2$ generated by $3$ forms of degree $d$.
The following conditions are equivalent:
\begin{itemize}
\item[{\rm (i)}] There exist two distinct minimal generating syzygies of degrees $r_1$ and $r_2$ such that
$r_1+r_2\geq d+1$.
\item[{\rm (ii)}] $I$ is not a perfect ideal.
\item[{\rm (iii)}] For any two distinct minimal generating syzygies of degrees $r_1$ and $r_2$, one has
$r_1+r_2\geq d+1$.
\end{itemize}
\end{lem}
\demo
(i) $\Rightarrow$ (ii)
This is clear: if $I$ is perfect then it is generated by the $2\times 2$ minors of a $3\times 2$ matrix
with columns of (standard) degrees that sum up to $d$.

(ii) $\Rightarrow$ (iii)
This is the core implication.
Suppose given two distinct minimal generating syzygies of degrees $r_1$ and $r_2$ such that
$r_1+r_2\leq d$.
We can assume that the given syzygies are columns of the minimal graded presentation matrix $\phi$ of $I$;
call $\psi$ the $3\times 2$ submatrix they form.
Dualizing into $R$, the kernel of the transposed $\phi^t$ is generated by a unique vector whose coordinates
give, on one hand, a set of minimal generators of $I$ and, on the other hand, are (up to sign) the $2\times 2$ minors
of the $2\times 3$ submatrix $\psi^t$ of $\phi^t$ divided by their $\gcd$.
It follows that this $\gcd$ must be unit, hence a set of generators of $I$ is the set of $2\times 2$ minors
of $\psi$.
This means that $I$ is a perfect ideal.

(iii) $\Rightarrow$ (i)
This is trivial.
\qed

\smallskip

Let ${\rm reg}(M)$ denote the Castelnuovo-Mumford regularity of a finitely generated graded module over
a standard polynomial ring over a field.

\begin{cor}
Let $R=\mathbb K[x,y,z]$ and $I\subset R$ denote an ideal of codimension $2$, minimally generated by $3$ forms of degree $d\geq 2$.
If ${\rm reg}(R/I)\leq \frac{3d-4}{2}$  then $I$ is perfect and equality ${\rm reg}(R/I)= \frac{3d-4}{2}$ holds, with $d$ necessarily even.
\end{cor}
\demo Pick two distinct minimal generating syzygies of $I$, of degrees $r_1$ and $r_2\geq r_1$.
Clearly, ${\rm reg}(R/I)\geq d+r_2-2\geq d+r_1-2$, hence $r_1+r_2\leq d$.
By Lemma~\ref{cod2_3gens}, $I$ is perfect.
But then ${\rm reg}(R/I)= d+r_2-2$ as it must be the case.
An immediate calculation now yields ${\rm reg}(R/I)\geq \frac{3d-4}{2}$.
It is obvious that $d$ must be even.
\qed

\smallskip

Suppose that $I$ is not perfect. The condition ${\rm indeg}(I^{\rm sat}/I)\geq d+1$ implies  the upper bound
${\rm reg}(R/I)\leq 2d-4$ (\cite[Theorem 1.5]{HS}).
This condition is necessary in general as the following
simple example of a gradient ideal shows.

\begin{exm} Let $I$ be the gradient ideal of $F=x(x^2+yz)$ (a conic and a transversal line).
We have $I=\langle 3x^2+yz,xz,xy\rangle$ and ${\rm reg}(R/I)=2$ which of course is not less or equal to $2d-4=0$.

Also $I^{\rm sat}=\langle x,yz\rangle$, and hence ${\rm reg}(R/I^{\rm sat})=1$.
Note that ${\rm indeg}(I^{\rm sat}/I)=1<d+1=3$ and that ${\rm st}(I)=2$.
\end{exm}
The above example has ${\rm reg}(R/I)=\frac{3d-4}{2}+\frac{{\rm st}(I)}{2}$.
In \cite[Corollary 1.7 (ii)]{HS} this value is shown to be attained whenever $I$ is non-perfect,
$I^{\rm sat}/I$ is generated in
a fixed degree and ${\rm indeg}(I^{\rm sat}/I)\geq d+1$.
It is remarkable that the last condition is very natural in the realm of gradient ideals
of reduced forms in $\mathbb K[x,y,z]$.
We wonder whether this bound is attained for gradient ideals under relaxed assumptions.
In fact, an answer is interesting enough even in the restricted case of gradients of line arrangements.

A couple of examples may illustrate this point.

\begin{exm}\label{arr1} Let $I\subset R$ denote the gradient ideal of the line arrangement $xyz(x+y+z)\in R$.
Then
$$I^{\rm sat}= \langle I,xyz\rangle =\langle x,y\rangle\cap\langle x,z\rangle\cap\langle y,z\rangle\cap
\langle x,y+z\rangle\cap\langle y,x+z\rangle\cap
\langle z, x+y\rangle$$
is minimally generated by $4$ forms of degree $3$ and one has st$(I)=1$.
The minimal free resolution of $R/I$ has the shape
$$0\rar R(-6)\lar R(-5)^3\lar R(-3)^3\lar R,$$
hence reg$(R/I)=3= \frac{3d-4}{2}+\frac{{\rm st}(I)}{2}$, with $d=3$.
\end{exm}

\begin{exm}\label{arr2}
Let $I\subset R$ denote the gradient ideal of the line arrangement
$$xyz(x+y)(x+z)(y+z)\in R.$$
Then
$$I^{\rm sat}= \langle I,g\rangle,$$
where $g$ is form of degree $6$.
Again,  st$(I)=1$, but this time around the additional two conditions stated in \cite[Corollary 1.7 (ii)]{HS}
hold true.
The minimal free resolution of $R/I$ has the shape
$$0\rar R(-9)\lar R(-8)^3\lar R(-5)^3\lar R,$$
hence reg$(R/I)=6= \frac{3d-4}{2}+\frac{{\rm st}(I)}{2}$, with $d=5$.
\end{exm}
 We note en passant that both examples are of ideals of {\em linear type}, to be introduced in the next
 subsection.

\subsection{Related algebras.}
To proceed, let $I=\langle f_1,f_2,f_3\rangle$ be a codimension 2 ideal in $R$. Then $I$ has a presentation
$$R^m\stackrel{\phi} \rightarrow R^3\rightarrow I\rightarrow 0.$$
Suppose that $\phi$ has columns the syzygies $(a_i,b_i,c_i),i=1,\ldots,m$.

\smallskip

Then the {\em symmetric algebra of $I$} is $${\mathcal S}_R(I)=R[T_1,T_2,T_3]/\langle a_1T_1+b_1T_2+c_1T_3,\ldots,a_mT_1+b_mT_2+c_mT_3\rangle.$$

Let $\mathcal R(I):=\oplus_{n=0}^{\infty}I^n$ denote the Rees algebra of $I$.
Then we have a surjective map of $R$-algebras $\psi:{\mathcal S}_R(I)\rightarrow \mathcal R(I)$ induced by
mapping $T_i\mapsto f_i$, where $f_i$ is viewed in degree $1$.
The algebra $\mathcal R(I)$ has a similar, but more complicated, presentation ideal.
Generators of these ideals are often referred to, respectively, as {\em equations of the symmetric algebra}
and {\em Rees equations} (or {\em equations of $I$}).

An ideal $I$ is called {\em of linear type} if the above map is injective (in other words, if ${\mathcal S}_R(I)\cong\mathcal R(I)$).
By \cite[Lemma 3.1 and Corollary 3.2]{AA} we have the following characterization (which is true for arbitrary number of variables):

\begin{prop}\label{old_thm} Let $I$ be an ideal of codimension 2 in $R=\mathbb K[x,y,z]$ minimally generated by three homogeneous forms.
Then the following are equivalent:
\begin{enumerate}
  \item $I$ is of linear type.
  \item $I$ localized at any of its minimal primes is a complete intersection.
  \item If $R^m\stackrel{\phi}\rightarrow R^3\rightarrow I\rightarrow 0$ is the minimal presentation of $I$,
  then the ideal generated by the entries of $\phi$ has codimension 3.
\end{enumerate}
\end{prop}

Actually, there is a more general statement (\cite[Proposition 3.7]{SV1}, also \cite[Remark 10.5]{HSV1} and \cite[Proposition 5.1]{SV2}):

\begin{prop} Let $I$ be an almost complete intersection in any Cohen-Macaulay ring $R$
such that ${\rm depth}\, R/I\geq \dim R/I-1$. If $I$ localized at any of its
minimal primes is a complete intersection {\rm (}i.e., generically a complete intersection{\rm )}, then ${\mathcal S}_R(I)\cong\mathcal R(I)$.
\end{prop}

The following result is very useful:

\begin{prop}\label{reg_seq}  {\rm \cite[Theorem 2.1]{t1}} Let $I\subset R$ be an ideal of codimension 2 minimally generated by three
homogeneous polynomials of the same degree. Then the following are equivalent:
\begin{enumerate}
  \item $I$ is perfect of linear type.
  \item $I$ has a syzygy forming a regular sequence of length $3$.
\end{enumerate}
\end{prop}

\smallskip

We note that Proposition~\ref{reg_seq} is not true in more than $3$ variables.
Similar restriction on the number of variables can be found in \cite[Theorem 1.7]{cs} which gives the
following syzygetical interpretation of locally complete intersections.

\begin{prop}\label{cox_schenck} Let $I\subset R=\mathbb K[x,y,z]$ be an ideal of codimension 2 minimally generated by three
homogeneous polynomials. Then the following are equivalent:
\begin{enumerate}
  \item $I$ is locally a complete intersection {\rm (}at its minimal primes{\rm )}
  \item The only syzygies  on $I$ with coordinates in $I^{\rm sat}$ are the Koszul syzygies.
\end{enumerate}
\end{prop}

Note that $I$ is perfect if and only if $I=I^{\rm sat}$. In the saturated case the above result has been established
in \cite[Proposition 2.3]{SV1} for higher codimension almost complete intersections.

Denoting $Z(I)\subset R^3$ (respectively, $K(I)$) the module of syzygies (respectively, Koszul syzygies) of $A,B,C$,
the second condition above reads as the equality
$Z(I)\cap I^{\rm sat}R^3=K(I).$
The relevance of the condition that $I$ be generically a complete intersection is discussed in
\cite[Corollary 2.6 and Question 2.7]{AHA} in its connection to the theorem of Dolgachev on reduced homaloidal plane
divisors.

Now, a basic condition discussed in (\ref{ssec2.1}) for an ideal generated by forms of the same degree $d$ is the
lower bound ${\rm indeg}(I^{\rm sat}/I)\geq d+1.$
This bound holds true for the gradient ideal a homaloidal polynomial
--  more generally, for any Cremona map of $\mathbb{P}^n$ with base ideal $I$.
For a proof of this fact see \cite[Proposition 1.2]{PanRusso} (also \cite[Proposition 2.5.2]{alberich} for the plane case).
It seems timely to pose the following
\begin{quest}\label{indeg}\rm
Let $I=J_F$ denote the gradient ideal of a reduced form $F\in \mathbb{K}[x,y,z]$ of degree $d+1$
satisfying the bound ${\rm indeg}(I^{\rm sat}/I)\geq d+1$.
When is $Z(I)\cap I^{\rm sat}R^3=K(I)$?
\end{quest}

\smallskip
If $I=I^{\rm sat}$ the bound is vacuously satisfied since ${\rm indeg}(\{0\})=+\infty$.
In this case, the question asks which saturated gradient ideals are {\em syzygetic} in the sense of \cite{MiRo}
(see also \cite[Section 2]{SV1}) -- i.e., such that the so-called syzygy part $Z(I)\cap IR^3/K(I)$
vanishes.
In Proposition~\ref{free_quintics} a many-parameter family of irreducible plane sextic curves with saturated gradient ideals
is given that fail for the above implication.
We now give a  $3$-parameter family of irreducible sextic curves, with non-saturated gradient ideals, for which
the implication fails in the strong sense that the respective gradient ideals are not even syzygetic.

\begin{prop}\label{indeg_example}\rm
Let $F=x^6+\alpha x^3y^3+\beta  x^2y^4 + \gamma y^5z\in R:=\mathbb K[x,y,z]$, where $\alpha,\beta,\gamma\in \mathbb{K}$ are nonzero.
Then
\begin{enumerate}
\item[{\rm (i)}] The gradient ideal $J_F\subset R$ is not saturated.
\item[{\rm (ii)}]  ${\rm indeg}(J_F^{\rm sat}/J_F)\geq 6$.
\item[{\rm (iii)}]  $Z(J_F)\cap J_F^{\rm sat}R^3\neq K(J_F)$ -- more strongly, $J_F$ is not syzygetic.
\end{enumerate}
\end{prop}
\demo Everything is computable with \cite{Macaulay}, by either: assuming the coefficients are indeterminates, then
virtually chopping off some superfluous relations  in the results pretending these indeterminates are scalars;
or else, taking random coefficients.
However, often a little preliminary hand calculation can tell us quite a bit in the right direction.

Set $I=J_F=\langle 6x^5+3\alpha x^2y^3+2\beta xy^4, \,3\alpha x^3y^2 + 4\beta x^2y^3+ 5\gamma y^4z,\, y^5\rangle$.

(i)
We claim:

\smallskip

$\bullet$ $(0,\, -y^3,\, 3\alpha x^3 + 4\beta x^2y+ 5\gamma y^2z)^t$ is a syzygy of $I$

\smallskip

$\bullet$ $x^2y^4\in I: \langle x,y,z\rangle\setminus I$.

\smallskip

The first is an immediate verification by looking at the shape of $F_y, F_z$.
For the second, the following three relations are directly confirmed:

{\small
\begin{eqnarray}\label{relations}\nonumber
x\cdot x^2y^4 & = & \underbrace{\frac{1}{3\,\alpha}\,y^2}_{v_1} F_y +\underbrace{\left(- \frac{1}{3\,\alpha\,\gamma}\,
(4\beta x^2+5\gamma\, yz)\right)}_{w_1} F_z=v_1F_y+w_1F_z\\
y \cdot x^2y^4  & = & \underbrace{\frac{1}{\gamma}\, x^2}_{w_2} F_z=w_2F_z\\ \nonumber
z \cdot x^2y^4  & = &\underbrace{\frac{-\alpha}{10\,\gamma}\,y^2}_{u_3} F_x+\underbrace{\left(\frac{1}{5\,\gamma}\, x^2-\frac{4\,\beta}{15\,\alpha\,\gamma}\,xy+\frac{16\,\beta^2}{45\,\alpha^2\,\gamma}\,y^2\right)}_{v_3} F_y\\ \nonumber
& + & \underbrace{\left(\frac{27\,\alpha^4-128\,\beta^3}{90\,\alpha^2\,\gamma^2}\,x^2+ \frac{\alpha\,\beta}{5\,\gamma^2}\,xy +\frac{4\,\beta}{\alpha\,\gamma}\,xz-\frac{16\,\beta^2}{9\,\alpha^2\,\gamma}yz\right)}_{w_3}F_z
=u_3F_x+v_3F_y+w_3F_z.\nonumber
\end{eqnarray}
}(The $u,v,w$'s introduced above will be revisited in Example~\ref{revisiting}).
The first two relations are visible, while the third one requires a calculation (it will be readily available with \cite{Macaulay}
by letting $\alpha,\beta,\gamma$ be honest indeterminates).

It remains to prove that $x^2y^4\notin I$. For this we resort to \cite{Macaulay}:
a Gr\"obner basis computation shows that $x^2y^4$ is not a multiple of any monomial in the initial ideal of $I$.
Thus, $I$ is not saturated.

(ii) We now show that the condition ${\rm indeg}(I^{\rm sat}/I)\geq 6$ holds.
Write $J:=\langle I,x^2y^4\rangle$. We have shown that $J\subset I:\langle x,y,z\rangle\subset I:\langle x,y,z\rangle ^{\infty}=I^{\rm sat}$.
It now suffices to show that $J$ is saturated, i.e., perfect.
But one can read (\ref{relations}) as syzygies of $J$; moreover, it is immediately checked that the
top $3\times 3$ minor of the matrix of these syzygies gives $\,y^2x^4$ (up to nonzero scalars).
From this it becomes apparent that $J$ is perfect with minimal resolution of shape
$$0\rar R(-7)^3\lar R(-5)^3\oplus R(-6)\lar J\rar 0.$$
It follows that $I^{\rm sat}=J$ and ${\rm indeg}(I^{\rm sat}/I)= 6$.

(iii)
To show the inequality $Z(I)\cap I^{\rm sat}R^3\neq K(I)$ we can apply Proposition~\ref{cox_schenck} by arguing
that $I$ is not a complete intersection locally at its unique minimal prime $(x,y)$; alternatively, one
can apply \cite{r}, or \cite[Section 1.3]{st} by showing that $F$ is not {\em weighted homogeneous} locally
at $(x,y)$ which is quite immediate here.

In addition, since $I$ is not saturated and since it is generated in degree $5$ it follows from Lemma~\ref{cod2_3gens}
that the initial standard degree of the syzygies is $3$ and the minimal number of syzygy generators is $3$.
From \cite[Proposition 1.6 (i)]{HS} the minimal free resolution of $I$ must be of the form
$$0\rar R(-9)\lar R(-8)^3\lar R(-5)^3 \lar I\rar 0.$$
(Note that this resolution has identical Betti numbers as the one in Example~\ref{arr2}.
However, $I$ is not of linear type as it was the case in loc. cit.)
Thus, we are in total control of the homological nature of $I$ and $I^{\rm sat}$.
To show the inequality $Z(I)\cap I^{\rm sat}R^3\neq K(I)$ in this variant it suffices to show that $I$ is not
a syzygetic ideal, i.e., the torsion $\delta(I)$ of the first Koszul homology module of $I$
does not vanish.
It is possible to show that $\delta(I)$ is cyclic generated by the residue of a syzygy of standard degree $7$,
hence yielding a Rees equation of $I$ of bidegree $(2,2)$ not coming from the equations of
the symmetric algebra of $I$.

\smallskip
A remnant problem is:

\begin{quest}\rm
Let $f\subset \mathbb{K}[x,y,z]$ denote a reduced Eulerian free divisor such that its gradient ideal $J_f$ is syzygetic.
When is $J_f$ of linear type?
\end{quest}
The question asks for such divisor for which the (Rees) equations not coming from the equations
of the symmetric algebra of $J_f$ have degree $\geq 3$.
At present we cannot pull out an example failing for this implication.

\smallskip

\subsection{When $I^{\rm sat}/I$ is the socle.} By definition, the socle of $R/I$ is $I:\langle x,y,x \rangle/I$.
Thus, $I^{\rm sat}/I$  coincides with the socle if and only  ${\rm st}(I)=1$.
Before we proceed discussing this interesting case, we analyze a special type of ideals.

\subsubsection{Ideals of Cramer type.} Let $R$ be a commutative ring, $n\geq 1$ an integer and  $\phi: R^n\rar R^n$ an $R$-homomorphism
of well-defined rank $n$  -- that is, $g:=\det(\phi)$ is a nonzerodivisor on $R$.
Given an element ${\bf u}\in R^n$ (target) there is a unique ${\bf t}\in R^n$ (source)
such that $\phi({\bf t})={\bf u}$.
Using matrix notation, letting $\psi:=[{\bf u}\,|\, \phi]: R^{n+1}\rar R^n$ denote the augmented map induced by ${\bf u}$,
by Cramer the coordinates of ${\bf t}$
in the canonical basis of the source, over the total quotient ring of $R$, are the $n$-minors of $\psi$ fixing the column ${\bf u}$
divided by $g$.

It is suggestive to call the ideal of $I\subset R$ generated by these minors $\{\Delta_1,\ldots, \Delta_n\}$ an ideal of {\em Cramer type
associated to the fixed column} $\bf u$.
We may call $g$ the {\em companion minor} of $I$.
Further, ``repeating'' a row of $\psi$ gives an $(n+1)\times (n+1)$ matrix with null determinant.
This implies that the rows of $\psi$ are syzygies of the Fitting ideal $I_n(\psi)=\langle I,g\rangle$.
If further $R$ is Noetherian and $I_n(\psi)$ has grade $\geq 2$, then one direction of the Hilbert--Burch theorem
says that
$$0\rar R^n \stackrel{\psi^t}{\lar} R^{n+1} \lar I_n(\psi)\rar 0$$
is a free resolution.

Here is one environment for Cramer ideals, based on a property of contents of ideals:

\begin{prop}\label{conductor_criterion}
Let $R$ be a Noetherian ring.
Let $I=\langle f_1,\ldots, f_n\rangle\subset R$ and $U=\langle u_1,\ldots, u_n\rangle\subset R$  be ideals, with
$I$ of grade $\geq 2$.
Let $g\in (I:U)$ and write
\begin{eqnarray}\label{condutor_eqs}
gu_1&=& a_{11}f_1+\cdots a_{1n}f_n\nonumber\\
&\vdots &\\
gu_n&=& a_{n1}f_1+\cdots a_{nn}f_n\nonumber
\end{eqnarray}
If the $n$-minors of the content matrix
\begin{equation}\label{matrix_g}
M_g:=\left(
\begin{array}{cccc}
u_1 &  a_{11}& \cdots &a_{1n}\\
\vdots &\vdots & \ddots & \vdots\\
u_n & a_{n1} & \cdots & a_{nn}
\end{array}
\right)
\end{equation}
generate an ideal of grade $\geq 2$, then $I$ is a Cramer ideal  associated to
the fixed column ${\bf u}:=(u_1,\ldots, u_n)^t$ with companion $g$.
Conversely, if $I$ is a Cramer ideal  associated to
the fixed column ${\bf u}:=(u_1,\ldots, u_n)^t$ with companion $g$ then $\langle I,g\rangle \subset I:U$
where $U$ is the ideal of $R$ generated by $\{u_1,\ldots, u_n\}$ and the $n$-minors of the resulting extended content matrix
generate an ideal of grade $\geq 2$.
\end{prop}
\demo
By an earlier observation, the ideal $I_n(M_g)$ has a free resolution
$$0\rar R^n \stackrel{M_g}{\lar} R^{n+1} \lar I_n(M_g)\rar 0.$$
Since $I_n(M_g)$ has grade at least $2$, dualizing into $R$ yields an exact sequence
$$0\rar R\stackrel{\rho}{\lar} {R^{n+1}}^* \stackrel{(M_g)^t}{\lar} {R^n}^*,$$
where the coordinates of a matrix of $\rho$ are the signed $n$-minors of $M_g$.
On the other hand, the equations (\ref{condutor_eqs}) show that $(M_g)^t$ is a matrix of syzygies of
$\{-g, f_1,\ldots, f_n\}$.
We can harmlessly assume that this matrix is a submatrix of a full syzygy $(n+1)\times p$ matrix $N$ of the latter elements.
Pick the corresponding  presentation of the ideal $J:=\langle -g, f_1,\ldots, f_n\rangle$ and dualize into $R$ to get an exact sequence
$$0\rar R\stackrel{\sigma}{\lar} {R^{n+1}}^* \stackrel{N^t}{\lar} {R^p}^*,$$
where the coordinates of $\sigma$ are the generators of $J$ since this ideal has grade $\geq 2$.
On the other hand, the coordinates of $\sigma$ are, up to signs and ordering, the $n$-minors of
any  $n\times (n+1)$ submatrix of $N^t$ of (maximal) rank $2$.
Since  $(M_g)^t$ sits as one of these, its $n$-minors must be, up to order and signs, coincide with
a set of generators of $J$.

This shows the first contention.
The converse statement is obvious.
\qed

\medskip

Ideals of Cramer type have been considered before in a larger realm (\cite{AnSi, Hochster}), but mainly in the generic
case -- say, when $R$ is a polynomial ring over a field and the entries of $\psi$ are independent indeterminates.
In this case,  the homological properties are quite known.
For example, by \cite[Theorem C]{AnSi} one has ${\mathcal S}_R(I)=\mathcal R(I)$, so $I$ is of linear type, and these algebras
are Cohen-Macaulay normal domains.

In the non-generic case, the linear type question doesn't have a precise answer. We'll see later in this subsection that the
gradient ideals considered in Example \ref{arr2} and Proposition~\ref{indeg_example} are Cramer type ideals, yet the first one is
of linear type and the second one is not.

Under extra conditions, \cite[Theorem D]{AnSi} gives a free resolution of $I$ in the non-generic case, based upon the well-known
Eagon-Northcott complex. Here is one version in the case where $\dim R=3$:

\begin{prop}\label{freeresCramer} Let $(R,\mathfrak{m})$ be a Cohen-Macaulay local ring, and let $I\subset \mathfrak{m}$ be as above.
Then the following are equivalent:
\begin{enumerate}
  \item The codimension of the ideal of $R$ generated by the coordinates of $\bf u$ is $3$ and $I_3(\psi)$ has codimension $2$.
    \item $I$ has a free resolution of the form $0\rightarrow R\rightarrow R^3\rightarrow R^3\rightarrow I\rightarrow 0$.
\end{enumerate}
\end{prop}

\medskip

We are now driven to apply the above general ideas in the case of an ideal $I=\langle f_1,f_2,f_3\rangle\subset R=\mathbb K[x,y,z]$ of
codimension $2$, minimally generated by three forms of degree $d\geq 2$.
Our basic assumption is that $I$ has saturation exponent $1$, i.e., that $I^{\rm sat}=I:\langle x,y,x \rangle$.
We harmlessly assume that  $gcd(f_i,f_j)=1, i\neq j$.

In this case, the  content matrix (\ref{matrix_g}) is of the form
\begin{equation}\label{sat_eqs}
M_g:=\left(
\begin{array}{cccc}
x&a_1 & b_1 & c_1\\
y&a_2& b_2 & c_2 \\
z&a_3& b_3 & c_3
\end{array}
\right),
\end{equation}
where the $a$'s,$b$'s and $c$'s  are forms of the same degree $s\geq 0$.
Fixing the set of generators of $I$, for given $g$ the matrix $M_g$ is not unique; however, any two such matrices
 are such that their respective transposes differ by a $3\times 3$ matrix of syzygies of the given set of generators of $I$.
All the same, we will say that the ideal $I$ satisfies the {\em generic saturation condition} (short: GSC) if there exists a
homogeneous element $g\in I^{\rm sat}\setminus I$ for which ${\rm codim}(I_3(M_g))=2$ (maximum possible).
We will also refer to such an element $g$ as a {\em saturation pivot}.

We collect additional properties of such ideals in the following

\begin{thm}\label{Cramer_vs_sat}
Let $I\subset R=\mathbb K[x,y,z]$ be an ideal of
codimension $2$, minimally generated by three forms of degree $d\geq 2$.
If $I$ has saturation exponent $1$ and satisfies condition {\rm GSC}, with saturation pivot $g$, then
\begin{enumerate}
\item[{\rm (i)}] $I$ is a Cramer ideal with fixed column $(x,y,z)^t$ and companion $g$ of degree $d-1+s$
for some $s\geq 1$.
\item[{\rm (ii)}] $I^{\rm sat}=\langle I,g\rangle$.
\item[{\rm (iii)}] $d$ is an odd integer.
\item[{\rm (iv)}] $\deg(g)=\frac{3}{2} (d-1)$ and ${\rm indeg}(I^{\rm sat}/I)\geq d+1$ if and only if $d\geq 5$.
\item[{\rm (v)}] The minimal graded free resolution of $R/I$ has the form
  $$0\rightarrow R\Bigl(-\frac{3}{2}(d+1)\Bigr)\rightarrow R^3\Bigl(-\frac{3d+1}{2}\Bigr)\rightarrow R^3(-d)\rightarrow R.$$
\end{enumerate}
Conversely, if $I$ is a non-saturated Cramer ideal with fixed column $(x,y,z)^t$ and companion $g$
then $I$ has saturation exponent $1$ and satisfies condition {\rm GSC}, with saturation pivot $g$.
\end{thm}
\demo
The assertion of (i) follows from Proposition~\ref{conductor_criterion}.
Here $s$ is the common degree of the entries of the matrix
\begin{equation}\label{g-matrix}
G:=\left(
\begin{array}{ccc}
a_1 & b_1 & c_1\\
a_2& b_2 & c_2 \\
a_3& b_3 & c_3
\end{array}
\right),
\end{equation}
which is $\geq 1$ since $g$ is not a scalar.

To prove (ii), let $J:=\langle I,g\rangle$.
By the proof of Proposition~\ref{conductor_criterion}, $J$ coincides with the ideal of $n$-minors of the
content matrix $M_g$.
Since the latter is saturated and contained in $I^{\rm sat}$, then it must be the case that they coincide.
This proves (ii).

To prove (iii), note that up to a nonzero scalar, $g$ coincides with the determinant of the matrix $G$
in (\ref{g-matrix}).
It follows that $\deg(g)=3s$, where $s$ is the degree of the entries of $G$.
On the other hand, we have seen that $\deg(g)=d-1+s$.
Therefore $d=2s+1$ is odd.

Then the value of $\deg(g)$ in (iv) follows immediately and it is easy to see that $\deg(g)\leq d$ if and only if $d\leq 3$.
Since $I^{\rm sat}=\langle I,g\rangle$ by (b), it is clear that ${\rm indeg}(I^{\rm sat}/I)=\deg(g)\geq d+1$
if and only if $d\geq 5$.

Finally, (v) follows from Proposition~\ref{freeresCramer} by reading the twists off the resolution in \cite{AnSi}
in the graded case.
Alternatively, one could proceed as in the proof of \cite[Theorem 2.14. (iii)]{HS}, working from the minimal graded resolution
of $R/\langle I,g\rangle$:
 $$0\rightarrow R^3\Bigl(-\frac{3d-1}{2}\Bigr)\stackrel{M_g}{\lar} R^3(-d)\oplus R\Bigl(-\frac{3}{2} (d-1)\Bigr)\rightarrow R.$$

\vskip .2in

For the converse statement, by Proposition~\ref{conductor_criterion} we have $g\in I:\langle x,y,x \rangle$ and the $n$-minors of
the extended content matrix $M_g$ generate an ideal of codimension $2$.

{\bf Claim}: $g\not\in I$.

To see this, suppose otherwise, letting $g=\alpha f_1+\beta f_2+\gamma f_3$. Plugging $g$ into the equations above
one obtains three syzygies of degree $s=\deg(g)+1-d=(d-1)/2$ of $f_1,f_2,f_3$:

$$\left(
\begin{array}{c}
a_1-x\alpha\\
b_1-x\beta\\
c_1-x\gamma
\end{array}
\right), \left(
\begin{array}{c}
a_2-y\alpha\\
b_2-y\beta\\
c_2-y\gamma
\end{array}
\right)\mbox{ and }\left(
\begin{array}{c}
a_3-z\alpha\\
b_3-z\beta\\
c_3-z\gamma
\end{array}
\right).$$
Since we are assuming that $I$ is not saturated and since $s+s=d-1$, from Lemma~\ref{cod2_3gens} we must conclude that
there is at most one minimal generating syzygy $(A,B,C)$ of degree $\leq d-1$, hence the three above syzygies  are multiple of this syzygy.
Letting $h_i, i=1,2,3$ denote the respective corresponding multiplying factors, one finds
$$g=\det \left(
\begin{array}{ccc}
h_1A+x\alpha & h_1B+x\beta & h_1C+x\gamma\\
h_2A+y\alpha & h_2B+y\beta & h_2C+y\gamma\\
h_3A+z\alpha & h_3B+z\beta & h_3C+z\gamma
\end{array}
\right)=0.$$

This contradicts the hypothesis $g\neq 0$ (included in the definition of a Cramer ideal).
This proves the claim.

Now, since $g\in I:\langle x,y,x \rangle\subset I^{\rm sat}$ then $\langle I,g\rangle \subset I^{\rm sat}$.
Since $\langle I,g\rangle=I_n(M_g)$ is unmixed, the equality $\langle I,g\rangle = I^{\rm sat}$ must be the case.
It follows that $I$ has satiety $1$ and satisfies $GSC$ with pivot $g$.
\qed

\medskip

\begin{exm}\label{revisiting} Let us look back at Proposition~\ref{indeg_example}, setting $f_1=F_x,f_2=F_x+F_y,f_3=F_x+F_z$ -- these are
minimal generators of $I$ such that $gcd(f_i,f_j)=1,i\neq j$. Then $I$ is a Cramer ideal with
associated fixed column $(x,y,z)^t$ and companion
$$
g:=-30\,\gamma^2\,\det
\left(
\begin{array}{ccc}
-v_1-w_1 & v_1 & w_1\\
-w_2& 0 & w_2 \\
u_3-v_3-w_3& v_3 & w_3
\end{array}
\right),
$$
where $u_i,v_j,w_k$ are the polynomials considered in Proposition~\ref{indeg_example}.
\end{exm}

\section{Homogeneous free divisors and the linear type property}

In the previous section we studied several homological properties of almost complete intersection ideals of codimension $2$
in $\mathbb K[x,y,z]$. As already suggested in that section, these results can be applied to the case when the ideal is
the gradient  ideal of a plane (projective or not) non-smooth divisor with defining polynomial being squarefree.

In this respect, the role of the symmetric  and  Rees algebras of the gradient ideal of an affine or homogeneous divisor establishes itself
naturally through the notion of {\em Koszul freeness}. Another natural carrier to introduce these algebras, in contrast to the classical
situation of completely reducible divisors, the linear type property of gradient ideals is not satisfied in general.

\subsection{Irreducible homogeneous divisors}

Our knowledge of irreducible divisors in $\mathbb P^n$ is at present very poor.
From what we could gather thus far, they seem to appear in isolated packets or special families.
In contrast to the case of divisors that, by definition, inherit natural combinatorial features
-- such as divisors which are union of lines (and conics, see \cite{st}) -- they leave not much choice but to resort to
the ideal theoretic and homological tool.

\subsubsection{Irreducible divisors in $\mathbb P^2$}
\medskip

We start with some elementary information on irreducible homogeneous divisors.

\subsubsection{Degree $\leq 4$.}

\medskip

\begin{prop}\label{free_divisors_degree4}
Let $\mathbb K$ denote an algebraically closed field.
An irreducible homogeneous  $f\in \mathbb K[x,y,z]$ of degree at most $3$ is a free divisor if and only if $\deg(f)=1$.
\end{prop}
\demo
If $\deg(f)=1$ with $f=ax+by+cz$ , say, $a\neq 0$, then clearly it is smooth.
An explicit determinantal form is
$$f=\det \left(
\begin{array}{ccc}
x & -b/\sqrt{a} & -c/\sqrt{a}\\
y& \sqrt{a} & 0 \\
z & 0 & \sqrt{a}
\end{array}
\right),
$$
hence $f$ is a free divisor -- note that it suffices to assume that $\mathbb K$ is closed under arbitrary square roots.
If $\deg (f)=2$ and $f$ is irreducible then ${\rm Proj}(\mathbb K[x,y,z]/(f))$ is smooth, hence the gradient ideal is a codimension $3$ complete
intersection (of linear forms).
Therefore, $f$ is not a free divisor.
If $\deg(f)=3$ and $f$ is irreducible, again either ${\rm Proj}(\mathbb K[x,y,z]/(f))$ is smooth, hence the gradient ideal is
a codimension $3$ complete
intersection (of quadrics), or else it has
a cusp or a node singularity. If $\fm=\langle x,y,x \rangle$ and $J_f$ is the corresponding gradient ideal, then
$xy\fm\subset J_f$ while $xy\not\in J_f$ in the cusp case,
and with a little more effort $xz\fm\subset J_f$ while $xz\not\in J_f$ in the node case.
Therefore in any case, ${\rm depth} R/J_f=0$ hence $J_f$ is not perfect.
\qed

\bigskip

Now assume that $\deg(f)=4$, with $f$ irreducible.
If $f$ is moreover rational then, based on \cite{Wall}, a classification of the corresponding families with fixed singular type
has been  given in \cite{AA} and it has been proved that the general member has a gradient ideal of linear type.
 Sweeping through the general member of each of these families with the help
of a computer program (Macaulay) one can show that the corresponding gradient ideal has homological dimension $2$
over $R=\mathbb K[x,y,z]$.

For irreducible quartics of higher genus a classification is roughly given in \cite{Nejad}.
Based on this classification, a second run with Macaulay shows that the corresponding gradient ideal has homological dimension $2$
over $R=\mathbb K[x,y,z]$.

Thus, it would seem that morally no irreducible homogeneous quartic is a free divisor.

\smallskip

In \cite{s3} an example was given of a family of irreducible homogeneous free divisors of degree $6$ in $\mathbb K[x,y,z]$
 -- a so-called {\em family of Cayley sextics}. This example is special in that quite some structure
 comes with the territory.

\subsubsection{Degree $\geq 5$.}

We now discuss a particular family of irreducible free divisors in $\mathbb{P}^2$ of any given degree $\geq 5$.

\begin{prop}\label{free_quintics} Let $F= y^{d-1}z+a_1x^d+a_2x^2y^{d-2}+a_3xy^{d-1}+a_4y^d\in \mathbb K[x,y,z]$,
with $d\geq 5$ and $a_1,a_2\neq 0$. Then:
\begin{itemize}
\item[{\rm (1)}] $V(F)$ is an irreducible free divisor for any $d\geq 5$.
\item[{\rm (2)}] $J_F$ is of linear type if and only if $d=5$.
\end{itemize}
\end{prop}
\demo First note that  any member of the family with $a_1\neq 0$ is an irreducible
polynomial in $\mathbb K[x,y,z]$. This is because on the affine piece $z=1$ the resulting polynomial
is the sum of two polynomials of successive degrees $d-1$ and $d$ with no common factor.

The partial derivatives of $F$ are:

\begin{eqnarray}
F_x&=&(d-1)a_1x^{d-1}+2a_2xy^{d-2}+a_3y^{d-1}\nonumber\\
F_y&=&\Bigl((d-2)a_2x^2+(d-1)(a_3x+z)y+da_4y^2\Bigr) y^{d-3}\nonumber\\
F_z&=&y^{d-1}.\nonumber
\end{eqnarray}

There is an obvious syzygy  $(0,\,-y^2,\,(d-2)a_2x^2+(d-1)(a_3x+z)y+da_4y^2)^t$.

Some slightly painful but straightforward calculation yields:

{\small
\begin{eqnarray*}
(d-2)a_2y^{d-3}F_x-da_1x^{d-3}F_y&=&p(x,y)F_z-d(d-1)a_1a_3x^{d-2}y^{d-2}-d(d-1)a_1x^{d-3}y^{d-2}z\\
a_3x^{d-2}y^{d-2}&=&\frac{a_3}{(d-2)a_2}x^{d-4}yF_y-q(x,y,z)y^{d-1}\\
x^{d-3}y^{d-2}z&=&\frac{1}{(d-2)a_2}x^{d-5}yzF_y-r(x,y,z)y^{d-1},
\end{eqnarray*}
} where $p(x,y)$ is a homogeneous polynomial of degree $d-3$ in the variables $x, y$, while $q(x,y,z)$ and $r(x,y,z)$ are homogeneous polynomials of degree $d-3$ in $\mathbb K[x,y,z]$, with $q(x,y,z)\in \langle x,y\rangle\mathbb K[x,y,z]$. Moreover, $x^{d-4}z$ (respectively, $x^{d-5}z^2$) is the only monomial in $q$ (respectively, $r$) for this power of $z$.

Since $y^{d-1}=F_z$, the above relations imply a syzygy of the shape
$$\Bigl((d-2)a_2y^{d-3},\, -da_1x^{d-3}+\cdots,\, \cdots+(d-1)x^{d-5}z^2\Bigr)^t$$ which is a syzygy of degree $d-3$.
Note that $2+(d-3)=d-1$. Therefore, to show that $V(F)$ is a free divisor it suffices, by Lemma~\ref{cod2_3gens}, to argue that these syzygies cannot be both multiple of one same minimal syzygy. But this is clear by the nature of the first coordinate of each of them.

Finally, again from the shape of the two generating syzygies, it is clear that their entries
generate an $\langle x,y,z\rangle-$primary ideal if and only if $d=5$.
Proposition~\ref{old_thm} then implies the second part of the statement.
\qed

\subsection{Weighted homogeneous divisors and their homogenization.} Let $f\in R:=\mathbb K[x_1,\ldots,x_n]$ be an arbitrary polynomial.
A natural question is to inquiry into what homological properties of the gradient ideal $J_f$ of $f$ are, so to say, inherited by
the gradient ideal $J_F$ of the homogenization $F\in R[t]$ of $f$ with respect to a new variable $t$.

\subsubsection{Recap of weighted homogeneous polynomials}

We will assume throughout that $f$ is not homogeneous in the standard grading of $R$, as otherwise the question is vacuous.

\begin{defn} A polynomial $f\in R:=\mathbb{K}[x_1,\ldots,x_n]$ is called {\em weighted homogeneous} {\rm (}or {\em quasi homogeneous}{\rm )}
if there exist positive rational numbers $0<a_1,\ldots,a_n \leq \frac{1}{2}$ {\rm (}called {\em rational weights}{\rm )} such that $$f(s^{a_1}x_1,\ldots,s^{a_n}x_n)=sf(x_1,\ldots,x_n)$$
\end{defn}
\noindent where $s$ is an indeterminate over $R$.

Taking the partial derivative with respect to $s$ and using the chain rule yields
$$\sum_{i=1}^n a_is^{a_i-1}x_if_{x_i}(s^{a_1}x_1,\ldots,s^{a_n}a_n)=f(x_1,\ldots,x_n).$$
Evaluating at $s=1$ gives
$$f=a_1x_1f_{x_1}+\cdots+a_nx_nf_{x_n}.$$

Conversely, let $f\in R$ be such that there exist
positive rational numbers $0<a_1,\ldots,a_n \leq \frac{1}{2}$ with $f=a_1x_1f_{x_1}+\cdots+a_nx_nf_{x_n}$.
Then by mapping $x_i\mapsto s^{a_i}x_i$, one obtains
$$f(s^{a_1}x_1,\ldots,s^{a_n}x_n)=s\,\frac{\partial}{\partial s}(f(s^{a_1}x_1,\ldots,s^{a_n}x_n)).$$
This is possible only when $$f(s^{a_1}x_1,\ldots,s^{a_n}x_n)=sC(x_1,\ldots,x_n),$$
for some $C(x_1,\ldots,x_n)\in R$.
Evaluating at $s=1$, yields $C(x_1,\ldots,x_n)=f(x_1,\ldots,x_n)$. So $f$ is a weighted homogeneous polynomial
of rational weights $a_1,\ldots,a_n$.

Instead of using rational weights it is common to use integer weights
$w_1,\ldots,w_n\geq 1$ with $\gcd=1$, in which case there exists a positive integer $\lambda$ such that
$f(s^{w_1}x_1,\ldots,s^{w_n}x_n)=s^{\lambda}f(x_1,\ldots,x_n)$ and $\lambda f=w_1x_1f_{x_1}+\cdots+w_nx_nf_{x_n}$.
The rational weights in the first form are then $a_i:=\frac{w_i}{\lambda},i=1,\ldots,n$.

\smallskip

The following is a key lemma in our analysis:

\begin{lem}\label{lin_syzygy} Let $f\in R$ be a weighted homogeneous polynomial of {\rm (}standard{\rm)} degree $d$ and
 rational weights $a_1,\ldots,a_n$. Let $F\in R[t]$ be the homogenization of $f$ with respect to a new variable $t$.
 Then the partial derivatives of $F$ admit the following linear syzygy
 $$(da_1-1)x_1F_{x_1}+\cdots+(da_n-1)x_nF_{x_n}-tF_t=0.$$
\end{lem}
\demo Dehomogenizing the Euler relation  $d\,F=x_1F_{x_1}+\cdots+x_nF_{x_n}+tF_t$
via $t\rightarrow 1$ gives
$$d\,f=x_1f_{x_1}+\cdots+x_nf_{x_n}+\widetilde{F_t},$$
where $\tilde{F_t}$ denotes the dehomogenization of $F_t$.

Using the assumption $f=a_1x_1f_{x_1}+\cdots+a_nx_nf_{x_n}$ yields
\begin{equation}\label{eq1}
(da_1-1)x_1f_{x_1}+\cdots+(da_n-1)x_nf_{x_n}-\widetilde{F_t}=0.
\end{equation}

On the other hand, writing $f=h_1+\cdots+h_d$, where $h_i\in R$ is homogeneous of degree $i$, one has
$F=h_1t^{d-1}+\cdots+h_{d-1}t+h_d$.
Then, (\ref{eq1}) implies the following relations
\begin{eqnarray}
(da_1-1)x_1(h_1)_{x_1}+\cdots+(da_n-1)x_n(h_1)_{x_n}&=&(d-1)h_1\nonumber\\
(da_1-1)x_1(h_2)_{x_1}+\cdots+(da_n-1)x_n(h_2)_{x_n}&=&(d-2)h_2\nonumber\\
&\vdots&\nonumber\\
(da_1-1)x_1(h_{d-1})_{x_1}+\cdots+(da_n-1)x_n(h_{d-1})_{x_n}&=&h_{d-1}\nonumber\\
(da_1-1)x_1(h_d)_{x_1}+\cdots+(da_n-1)x_n(h_d)_{x_n}&=&0.\nonumber
\end{eqnarray}

Multiplying each of these equations by a suitable power of $t$ and adding up the resulting relations
yields
$$(da_1-1)x_1F_{x_1}+\cdots+(da_n-1)x_nF_{x_n}-tF_t=0,$$
as was to be shown.
\qed

\medskip

\subsubsection{Homogenization of binary weighted homogeneous polynomials}
We have the following immediate characterization of weighted homogenous polynomials in two variables, with positive integral weights.

\begin{prop}\label{binary_quasi_homog}
The weighted homogeneous polynomials in $\mathbb K[x,y]$ with positive integral weights form a family with typical member of the form
$$ f=c_x x^{sq}+c_y y^{sp}+\sum_{1\leq r<s} c_r x^{(s-r)q}y^{rp},$$
for a triple $(p,q,s)$, with $p,q,s\in {\mathbb N}_{>0}$,  such that $\gcd(p,q)=1$
and, moreover, weights $pe,qe$, with $e\geq 1$ arbitrarily chosen.
\end{prop}
\demo
The proof follows immediately from the definitions.
\qed

\medskip

Furthermore, note the following additional attributes of a binary weighted homogeneous polynomials:

$\bullet$ Homogeneous: case $p=q=1$

$\bullet$ Irreducibility: at least both $c_x\neq 0, c_y\neq 0$

$\bullet$ Partial irreducibility: neither $x$ nor $y$ is a factor of $f$ if and only if at least one of the
rational weights is $1/\deg(f)$ (standard degree)

$\bullet$ Smoothness at $(0,0)$: irreducible plus $sq=1$ or $sp=1$ (i.e., $s=1$ and either $p=1$ or $q=1$).

\medskip

For the case of divisors on $\mathbb P^2$ obtained by homogenization of weighted homogeneous polynomials in two variables, we obtain
a complete homological characterization as follows.

\begin{thm}\label{three_vars} Let $f\in R:=\mathbb K[x,y]$ be a {\rm (}non-homogeneous{\rm )} weighted homogeneous
polynomial of degree $d\geq 3$, with $f=axf_x+byf_y, a,b>0$. Let $F\in S:=R[z]$ be the homogenization of $f$ with respect
to a new variable $z$. Then
\begin{enumerate}
\item If $a\neq\frac{1}{d}$ and $b\neq\frac{1}{d}$, then $J_F$ is perfect of linear type with graded minimal free resolution:
$$0\rightarrow S(-(2d-3))\oplus S(-d)\rightarrow S^3(-(d-1))\rightarrow J_F\rightarrow 0.$$
\item If $a=\frac{1}{d}$ or $b=\frac{1}{d}$, then $J_F$ is an ideal of linear type, but it is not perfect. In this case,
$J_F$ has a graded minimal free resolution:
$$0\rightarrow S(-(2d-1)){\rightarrow} S(-d)\oplus S(-2(d-1))^2 {\rightarrow}S(-(d-1))^3{\rightarrow} J_F\rightarrow 0.$$
\end{enumerate}
\end{thm}
\demo (1) By Lemma \ref{lin_syzygy},  $((da-1)x,(db-1)y,-z)$ is a syzygy on $J_F$.
 Clearly, its coordinates form a regular sequence of length $3$. Therefore, the desired result follows from Proposition~\ref{reg_seq}.
\vskip .2in
(2) Since the data are symmetric, we may assume that $a=1/d$ (then $b\neq 1/d$, as we are assuming that $f$
is not homogeneous
 in the standard grading of $R$). Then $x^{d-1}$ occurs as a pure power in $F_x$. Also, Lemma \ref{lin_syzygy} gives
 $$(db-1)yF_y=zF_z.$$
 Denote $\mathbf{z}:=(0,(db-1)y, -z)^t$ the corresponding linear syzygy.

$I_1(\phi)$ contains $y,z$, where $\phi$ denotes a minimal graded presentation matrix of $J_F$. On the other hand, $F_x$ contains a
nonzero pure power in $x$, while its remaining terms belong to the ideal $\langle y,z\rangle$. Because of any of the Koszul syzygies
involving $F_x$ this forces $I_1(\phi)$ to be $\langle x,y,z\rangle$-primary. But, since $J_F$ is an almost complete intersection,
by Proposition~\ref{old_thm}, this implies that $J_F$ is of linear type.
\medskip

To prove the second part we argue in the following way. If $F$ is free then we have that $\mathbf{z}=(0,(db-1)y, -z)^t$ and
$(A,B,C)^t$ basis for the syzygies module of $J_F$, for some syzygy $(A,B,C)$: as $\deg(F)\geq 3$, we have that $\mathbf{z}$ is
the only linear syzygy so it must be a part of the basis.

Since $(db-1)yF_y=zF_z$, then $$F_y=zH, F_z=(db-1)yH,$$ for some $H\in S$. The $2\times 2$ minors of the matrix of the two syzygies
will give linear combinations of $F_x, F_y, F_z$. So $F_x \in \langle y,z\rangle$. Contradiction.

\medskip

Consider the two Koszul syzygies $\mathbf{z}_{12}:=(F_y,-F_x,0)^t$ and $\mathbf{z}_{13}:=(F_z,0,-F_x)^t$.
Clearly, $\mathbf{z}, \mathbf{z}_{12}, \mathbf{z}_{13}$ form a minimal set of generators of the $S$-submodule
of $S^3$ that they generate.

A direct calculation shows that the vector $({F}_x, -z, (db-1)y)^t$ is a syzygy of the set
$\{\mathbf{z}, \mathbf{z}_{12}, \mathbf{z}_{13}\}$.
This way one has a complex of $R$-modules
$$0\rightarrow S(-(2d-1)) \stackrel{\psi_2}{\rightarrow} S(-d)\oplus S(-2(d-1))^2 \stackrel{\psi_1}{\rightarrow}
S(-(d-1))^3 \stackrel{\psi_0}{\rightarrow} S,$$
where Im$(\psi_2)=(F_x, \beta z, \alpha y)^t$, Im$(\psi_1)=S\mathbf{z}+S \mathbf{z}_{12}+S \mathbf{z}_{13}$
and Im$(F_x,F_y,F_z)$.

To see that it is a resolution of its cokernel $J_F$, we apply the well-known acyclicity criterion as formulated
in \cite[Theorem 20.9]{E}. The ranks of the maps clearly add up right. Thus, it suffices to check the heights of the
appropriate Fitting ideals.

Since $F_x$ contains a nonzero pure power in $x$, while its remaining terms belong to the ideal $\langle y,z\rangle$, then the
height of $I_1(\psi_2)$ is $3$.

As to $\psi_1$, it has the $2$-minor $F_x^{\,2}$ which has a non-vanishing pure power term in $x$, while it has another
$2$-minor belonging to the ideal $\langle y,z\rangle$. The two cannot have a proper common factor, hence $I_2(\psi_1)$
has height at least $2$.
\qed

\medskip

It may be of some convenience to state the theorem in terms of free divisors
and the partial reducibility attribute mentioned earlier:

\begin{cor}
Let $f\in R:=\mathbb K[x,y]$ stand for a  weighted homogeneous
polynomial of degree $d\geq 3$, with positive rational weights $a,b$ and let $F\in R[z]$ denote its homogenization.
Then $F$ is a free divisor if and only if either $f$ is {\rm (}standard{\rm )} homogenous or else $f$ is
divisible by $x$ or by $y$.
\end{cor}

\subsubsection{Coneing.}
From the previous section we have seen that the homogenization of a weighted homogeneous polynomial in two variables leads
to divisors of linear type, but most of the time not free (e.g., if $f$ is irreducible). So freeness is not preserved
under homogenization in this low dimension.

Actually, by and large, homogenization is perhaps not the right operation to look at. An alternative operation,
enjoying good properties at least for divisors which are hyperplane arrangements,  is {\em coneing}.
Let $f\in R:=\mathbb K[x_1,\ldots,x_n]$ be a polynomial, which is not homogeneous. Let $F\in R[t]$ be
the homogenization of $f$ with respect to a new variable $t$. Then {\em the cone over $f$} is the homogeneous
polynomial ${\rm cone}(f):=tF\in R[t]$.

As a consequence to Lemma~\ref{lin_syzygy}, we have the following result.

\begin{prop}\label{free_cones} Let $f\in R$ be a weighted homogeneous polynomial of degree $d$ and with rational weights $a_1,\ldots,a_n$.
Denote by $G:={\rm cone}(f)\in R[t]$, the cone over $f$. Let $F=\frac{G}{t}$. If the ideal $I_F=\langle F_{x_1},\ldots,F_{x_n}\rangle$
is a perfect ideal of codimension $2$, then the gradient ideal $J_G$ of $G$ is perfect.
\end{prop}
\demo Let $F$ be the homogenization of $f$ with respect to the variable $t$. Then $G=tF$. The partial derivatives of $G$ are
$$G_{x_1}=tF_{x_1}, \ldots, G_{x_n}=tF_{x_n}, G_t=F+F_t.$$

Because $dF=x_1F_{x_1}+\cdots+x_nF_{x_n}+tF_t$, and from the linear syzygy obtained in Lemma \ref{lin_syzygy}, one gets that $$G_t=(da_1+a_1-1)x_1F_{x_1}+\cdots+(da_n+a_n-1)x_nF_{x_n}.$$

$I_F$ is a perfect ideal of codimension 2, so let

$$\left(
\begin{array}{ccc}
h_{1,1} & \cdots & h_{1,n-1}\\
 & \ddots &\\
h_{n,1} & \cdots & h_{n,n-1}
\end{array}
\right)
$$

\noindent be the $n\times (n-1)$ matrix whose (signed, ordered) $(n-1)$-minors are $\{F_{x_1}, \ldots,  F_{x_n}\}$.
Then the partial derivatives of $G$ are the (signed) $n$-minors of the $(n+1)\times n$ matrix

$$\left(
\begin{array}{ccc@{\quad\vrule\quad}c}
\raise5pt\hbox{$h_{1,1}$}&\raise5pt\hbox{$\cdots$}&\raise5pt\hbox{$h_{1,n-1}$}&\raise5pt\hbox{$(da_1+a_1-1)x_1$}\\[-6pt]
&\raise5pt\hbox{$\ddots$}&&\raise5pt\hbox{$\vdots$}\\[-6pt]
\raise5pt\hbox{$h_{n,1}$}&\raise5pt\hbox{$\cdots$}&\raise5pt\hbox{$h_{n,n-1}$}&\raise5pt\hbox{$(da_n+a_n-1)x_n$}\\
\multispan4\hrulefill\\[0.1pt]
0&\cdots&0& t
\end{array}
\right)
$$
\qed

\medskip

\begin{quest} If $f\in R$ is a weighted homogeneous polynomial such that $V(f)$ is a free non-smooth reduced divisor, is it
true that the ideal $I_F$ occurring in Proposition \ref{free_cones} is perfect of codimension 2?
\end{quest}

If $f$ is in two variables, the answer is yes, as we can see immediately at the beginning of the proof of the next result.

\begin{cor} Let $f\in R:=\mathbb K[x,y]$ be a weighted homogeneous polynomial of degree $d\geq 3$, with rational weights $a,b$.
Suppose $V(f)\subset \mathbb A^2$ is a reduced non-smooth divisor. Denoting $G={\rm cone}(f)\in S:=R[z]$, then $J_G$ is a
codimension 2 perfect ideal of linear type with graded minimal free resolution
$$0\rightarrow S(-(d+1))\oplus S(-(2d-1))\rightarrow S(-d)^3\rightarrow J_G\rightarrow 0.$$
\end{cor}
\demo Let $F=\frac{G}{z}$ be the homogenization of $f$. From Proposition \ref{free_cones}, it is enough to prove that
${\rm codim}(\langle F_x,F_y\rangle)=2$ (so $I_F$ is a complete intersection), which is equivalent to show that
$F_x$ and $F_y$ do not have a common factor.

If $F_x$ and $F_y$ have a common factor that is not a pure power of $z$, then by dehomogenization, $f_x$ and $f_y$
will have a non-constant common factor, contradicting ${\rm codim}(\langle f_x,f_y\rangle)=2$.

If $F_x$ and $F_y$ have a common factor that is a pure power of $z$, say $z^u, u\geq 1$, then by dehomogenizing,
$\deg(f_x)\leq d-1-u$ and $\deg(f_y)\leq d-1-u$, contradicting the fact that $\deg(f)=d$, and hence at least one
of its partial derivatives must have degree equal to $d-1$.

\medskip

We have that $$G_x=zF_x,G_y=zF_y,G_z=F+zF_z.$$ Using the Euler relation $dF=xF_x+yF_y+zF_z$, and the syzygy
$(da-1)xF_x+(db-1)yF_y=zF_z$ from Lemma \ref{lin_syzygy}, we obtain
$$G_z=(d(a+1)-1)xF_x+(d(b+1)-1)yF_y.$$
Multiplying this by $z$, and since $d(a+1)-1\neq 0$ and $d(b+1)-1\neq 0$ (otherwise $a<0$ or $b<0$), we obtain
the syzygy
$$(d(a+1)-1)x,(d(b+1)-1)y,-z),$$
of $G_x,G_y,G_z$.
The entries in this syzygy form a regular sequence of length 3, and therefore the conclusion of the theorem is shown.
\qed

\subsection{New and old examples}

\subsubsection{Adding free homogeneous divisors to irreducible homogeneous divisors}

A natural question is whether given any irreducible homogeneous polynomial
$g\in R:=\mathbb K[x_1,\ldots,x_n]\, (n\geq 3)$ there exists a
suitable reduced homogeneous free divisor $h\in R$ such that $f:=gh$ is a free divisor.
Next question would ask what is the minimal number of irreducible components such an $h$ should have.
Of course, in the best of the cases one would hope to accomplish this with $h$ itself irreducible.

The following propositions give a sense of how this addition problem works.
We assume that $n\geq 3$ throughout this part.

\begin{prop}\label{addition}
Let $g=x_1^{r_1}x_2^{r_2}\cdots x_{n-1}^{r_{n-1}}-x_n^d\in R$, with $r_1+\cdots + r_{n-1}=d$ and $r_i\neq 0$ for $i=2,\ldots,n-1$,
and let $h:=x_{i_1}\cdots x_{i_{n-2}}$, where $\{i_1,\ldots,i_{n-2}\}\subset \{1,\ldots, n-1\}$.
Then
\begin{enumerate}
\item[{\rm (i)}] $g$ is a nonfree divisor.
\item[{\rm (ii)}] $f:=gh$ is a free divisor of linear type.
\end{enumerate}
\end{prop}
\demo
(i)
To prove that $g$ is not a free divisor, we use the following general obvious

{\sc Claim.} Let $I\subset R$ be a codimension $2$ ideal generated by at least three independent forms of the same degree $\geq 2$.
If some Koszul relation among these forms is a minimal syzygy generator then $I$ is not perfect.

To apply in our case, note that by the shape of $g$, the partial derivative $g_{x_n}$ is a pure power of $x_n$ while
none of the others  involves $x_n$. Then $x_n$ is a nonzerodivisor modulo the latter, hence at least one Koszul
relation involving $x_n$ is a minimal syzygy generator.

(ii)
Without loss of generality we assume that $\{i_1,\ldots,i_{n-2}\}=\{1,\ldots, n-2\}$.
One then has
$$f=x_1^{r_1+1}x_2^{r_2+1}\cdots x_{n-2}^{r_{n-2}+1}x_{n-1}^{r_{n-1}}- x_1x_2\cdots x_{n-2}x_n^d,$$
from which we see that
\begin{eqnarray*}
f_{x_1}= (r_1+1)x_1^{r_1}x_2^{r_2+1}&\cdots &x_{n-2}^{r_{n-2}+1}x_{n-1}^{r_{n-1}}- x_2\cdots x_{n-2}x_n^d  \\
f_{x_2}=  (r_2+1)x_1^{r_1+1}x_2^{r_2}x_3^{r_3+1}&\cdots & x_{n-2}^{r_{n-2}+1}x_{n-1}^{r_{n-1}} - x_1x_3\cdots x_{n-2}x_n^d  \\
&\vdots &\\
f_{x_{n-2}}= (r_{n-2}+1)x_1^{r_1+1}x_2^{r_2+1}&\cdots & x_{n-3}^{r_{n-3}+1}x_{n-2}^{r_{n-2}}x_{n-1}^{r_{n-1}} - x_1x_2\cdots x_{n-3}x_n^d  \\
f_{x_{n-1}}=r_{n-1}\,x_1^{r_1+1}x_2^{r_2+1}&\cdots & x_{n-3}^{r_{n-3}+1}x_{n-2}^{r_{n-2}+1}x_{n-1}^{r_{n-1}-1}\\
f_{x_{n}}=  d\,x_1x_2\cdots x_{n-2}x_n^{d-1}
\end{eqnarray*}
From the last two partial derivatives, which are monomials, we get a {\em reduced} Koszul relation of degree $d-1$.
Suppose we have found $n-2$ additional minimal syzygies which are linear.
Then we argue as in the proof of Lemma~\ref{cod2_3gens} ((ii) $\Rightarrow$ (iii)) by dualizing this partial matrix into $R$
so that the $(n-1)$-minors of the transpose divided by their $\gcd$ gives the partial derivatives.
But since these minors have degree $d-1+n-2=d+n-3=\deg(f)-1$, we conclude that the syzygies obtained so far
already generate all syzygies, hence the gradient ideal is a (codimension $2$) perfect ideal.

Now, for any $i\in\{1,\ldots, n-2\}$, the following is a linear syzygy
$$\Bigl(0,\ldots,0, \underbrace{x_i}_{i\mbox{\rm th slot}},0, \ldots,0, -\frac{r_i+1}{r_{n-1}}\,x_{n-1}, -\frac{1}{d}\,x_n\Bigr)^t,$$
as one verifies straightforwardly.
Because of the echelon positioning of the $i$th coordinate, these $n-2$ syzygies are independent (hence, minimal) generators.

To complete the proof, we argue about the linear type property.
Since $I$ is a codimension $2$ perfect ideal, the Koszul homology modules of $I$ are Cohen--Macaulay (\cite{AvHe}).
Therefore, by \cite[Theorem 9.1]{HSV1} it suffices to verify that $\mu(I_P)\leq \dim R_P$ for every prime ideal
$P\supset I$. This condition is equivalent to checking the so-called $(F_1)$ property of the Fitting ideals of the
syzygy matrix of $I$:

$$\phi:=\left(
\begin{array}{cccc@{\quad\vrule\quad}c}
\raise5pt\hbox{$x_1$}&&&&\\[3pt]
&\raise5pt\hbox{$x_2$}&&&\\
&&\raise5pt\hbox{$\ddots$}&&\\
&&&\raise5pt\hbox{$x_{n-2}$}&\\
\multispan5\hrulefill\\[0.1pt]
c_1x_{n-1} & c_2x_{n-1} & \cdots & c_{n-2}x_{n-1}& cx_n^{d-1}\\[2pt]
c'_1x_{n-1} & c'_2x_{n-1} & \cdots & c'_{n-2}x_{n-1}& \mathfrak{m}
\end{array}
\right),
$$
\newline
where the empty slots have null entries, $c_j,c'_j, c$ are suitable nonzero scalars and $\mathfrak{m}$ denotes a suitable monomial.
The property to be verified is that, for every $1\leq t\leq n-1$, the codimension of the ideal $I_t(\phi)$
is at least $n+1-t$.

Fix $1\leq t\leq n-1$.
Considering the main diagonal of the $(n-2)\times (n-2)$ upper left matrix, extended by the entry $cx_n^{d-1}$,
 leads us to look at the ideal generated by the $t$-products of the regular sequence $\{x_1,\ldots,x_{n-2}, cx_n^{d-1}\}$,
which is well-known to have codimension $n-1-t+1=n-t$.
Thus, we only need an additional $t$-minor effectively involving the remaining variable $x_{n-1}$.
But this is easily obtained using one of the last two rows of $\phi$.
\qed

\medskip

A sort of generalization of Proposition~\ref{addition} is the following:

\begin{thm}\label{addition2} Let $g=x_n^d-x_{n-1}^mh\in \mathbb K[x_1,\ldots,x_n], d>m$, with $h\in\mathbb K[x_1,\ldots,x_{n-2}]$
a reduced homogeneous polynomial of degree $d-m\geq 2$.
\begin{enumerate}
\item[{\rm (i)}] $g$ is a nonfree divisor.
\item[{\rm (ii)}] If $h$ is a free divisor then $f:=gh$ is a free divisor {\rm ;} if, moreover, $J_h$ is an ideal
of linear type in $\mathbb K[x_1,\ldots,x_{n-2}]$ then $J_f$ is of linear type.
\end{enumerate}
\end{thm}
\demo (i) Same argument as in proof of Proposition~\ref{addition} (i) gives a minimal generating Koszul syzygy on the
partial derivatives of $g$. From (the generalized version) Lemma~\ref{cod2_3gens}, obviously $g$ cannot be free.

\medskip

(ii) By assumption, $f=x_n^dh-x_{n-1}^mh^2$.
A simple calculation  gives
$$f=\frac{1}{d}x_nf_{x_n}+\frac{1}{m}x_{n-1}f_{x_{n-1}}.$$
Using the Euler's formula on $f$ we obtain a linear syzygy
$${\rm{\bf s}_0}=(x_1,\ldots,x_{n-2},(1-\frac{2d-m}{m})x_{n-1},(1-\frac{2d-m}{d})x_{n})^t.$$

Since $h$ is free there exists a basis of $n-3$ syzygies on $J_h$ with degrees adding up to $d-m-1$.
On the other hand,
\begin{equation}\label{H_factor}
f_{x_i}=H\,h_{x_i}, \,\, H:=(x_n^d-2x_{n-1}^mh),i=1,\ldots,n-2,
\end{equation}
as one easily checks.
Therefore, these syzygies extend to syzygies on $J_f$ by adding two more entries equal to zero.
Let ${\rm{\bf s}_1},\ldots,{\rm{\bf s}_{n-3}}$ denote these syzygies.

Finally, as in the previous argument, since $f_{x_{n-1}}=-mx_{n-1}^{m-1}h^2, f_{x_n}=d\,x_n^{d-1}h$ admit
the common factor $h$,
there is a syzygy of degree $d-1$:
 $${\rm{\bf s}_{n-2}}=(0,\ldots,0,dx_n^{d-1},mx_{n-1}^{m-1}h)^t.$$
So far, we have the following partial matrix of syzygies of $J_f$:

\begin{equation}\label{matrix_of_Jf}
\phi:=\left(
\begin{array}{ccc@{\quad\vrule\quad}cc}
&& & x_1 & 0\\[3pt]
 & \raise5pt\hbox{$\phi_h$} & & \raise5pt\hbox{$\vdots$} &\raise5pt\hbox{$\vdots$}\\
&& & \raise5pt\hbox{$x_{n-2}$} & \raise5pt\hbox{$0$}\\
\multispan5\hrulefill\\[0.5pt]
0&\cdots & 0 & -2\frac{(d-m)}{m} x_{n-1} & d\,x_n^{d-1}\\[5pt]
0&\cdots & 0 & -\frac{(d-m)}{d}x_{n} &mx_{n-1}^{m-1}h
\end{array}
\right),
\end{equation}
where $\phi_h$ denotes the syzygy matrix of $J_h$.
We claim that $\phi=\phi_f$. For it we use the same argument as before by embedding $\phi$ into $\phi_f$
and taking the transpose $\phi_f^t$. By the shape of the syzygies, $\phi^t$ has rank $n-1$, hence its $(n-1)$-minors
divided by their $\gcd$ give the elements of $J_f$.
But the degrees of ${\rm{\bf s}_0},{\rm{\bf s}_1},\ldots,{\rm{\bf s}_{n-3}},{\rm{\bf s}_{n-2}}$ add up to $1+(d-m-1)+(d-1)=2d-m-1$
which is the degree of the generators of $J_f$
Therefore, the $\gcd$ is $1$, thus implying the minors already generate $J_f$.

\smallskip

To prove the assertion about the linear type  property, we now inspect the syzygy matrix $\phi_f$  of $J_f$,
which we now know is given by (\ref{matrix_of_Jf}).
As explained earlier it suffices to show that for every $1\leq t\leq n-1$, the codimension of the ideal $I_t(\phi_f)$
is at least $n-t+1$.
But, this condition is also necessary, hence applying it to $\phi_h$ we get the inequalities
\begin{equation}\label{F1_of_h}
{\rm codim}(I_t(\phi_h))\geq n-3-t+2=n-t-1, \quad i=1,\ldots, n-3,
\end{equation}
as ideals in the ring $\mathbb K[x_1,\ldots,x_{n-2}]$.

Consider the following $(n-1)\times (n-2)$ submatrix of $\phi_f$
$$\psi_1=\left(
\begin{array}{ccc@{\quad\vrule\quad}c}
&& & x_1 \\[3pt]
 & \raise5pt\hbox{$\phi_h$} & & \raise5pt\hbox{$\vdots$}\\
&& & \raise5pt\hbox{$x_{n-2}$}\\
\multispan4\hrulefill\\[0.1pt]
0&\cdots & 0 & cx_{n-1} \\[2pt]
\end{array}
\right).
$$
Note that the upper $(n-2)$-minor of $\psi_1$ is just $h$ (up to a nonzero scalar multiple) since the corresponding
submatrix is the Saito matrix of $h$.
Since, on the other hand, $c\neq 0$, one has $I_{n-2}(\psi_1)=\langle x_{n-1}J_h,  h\rangle$.
Now any prime ideal containing $I_{n-2}(\psi_1)$ must either contain $\langle x_{n-1},h\rangle$ or else
$\langle J_h,h\rangle =J_h$. Since $J_h\subset k[x_1,\ldots,x_{n-2}]$, $x_{n-1}$ is a nonzerodivisor on $J_h$,
and since $h$ is assumed to be reduced, $J_h$ has codimension $\geq 2$ (actually $=2$).
Therefore, $I_{n-2}(\psi_1)$ has codimension $2$.

We can extend this result to any  $1\leq t \leq n-2$ as follows:
\begin{equation}\label{psi1}
{\rm codim}(I_t(\psi_1))\geq n-t
\end{equation}
in $\mathbb K[x_1,\ldots,x_{n-2},x_{n-1}]$.

To see this, let $\mathfrak{p}$ be a minimal prime over $I_t(\psi_1)$. We have $x_{n-1}I_{t-1}(\phi_h)\subset I_t(\psi_1)\subset \mathfrak{p}$.
If $x_{n-1}\notin \mathfrak{p}$, then $I_{t-1}(\phi_h)\subset \mathfrak{p}$, and therefore
${\rm codim}(\mathfrak{p})\geq (n-2)+1-(t-1)=n-t$. If $x_{n-1}\in \mathfrak{p}$, since $I_t(\phi_h)\subset I_t(\psi_1)$, then $\langle x_{n-1},I_t(\phi_h)\rangle\subset \mathfrak{p}$. But $I_t(\phi_h)$ is an ideal in $\mathbb K[x_1,\ldots,x_{n-2}]$,
of codimension $\geq n-t-1$. So ${\rm codim}(\mathfrak{p})\geq n-t$ in this case as well.
This proves the claim.

Now, for $1\leq t\leq n-1$,  consider the subideal of $I_t(\phi_f)$ generated by the $t$-minors involving the
lower right $2\times 2$ submatrix of $\phi_f$:
\begin{equation*}
\left(
\begin{array}{cc}
 -2\frac{(d-m)}{m} x_{n-1} & d\,x_n^{d-1}\\[5pt]
 -\frac{(d-m)}{d}x_{n} &mx_{n-1}^{m-1}h
\end{array}
\right)
\end{equation*}
Note that its determinant is the factor $H$ in (\ref{H_factor}).
Thus, for $t=n-1$, we are recovering the subideal as in (\ref{H_factor}).
In general, this subideal is precisely $H\, I_{t-2}(\phi_h)$.
One has for $1\leq t\leq n-2$:
$$ I_t(\psi_1)+H\, I_{t-2}(\phi_h)\subset I_t(\phi_f).$$
Let $\mathfrak{p}$ be prime containing $I_t(\phi_f)$.
If $H\in \mathfrak{p}$ then $I_t(\psi_1) +\langle H\rangle\subset \mathfrak{p}$; but $H$ is a nonzerodivisor modulo $I_t(\psi_1)$ since
it has a term which is a pure power in $x_n$.
Therefore, the codimension of $\mathfrak{p}$ is $\geq n-t+1$, where we have used the estimate in (\ref{psi1}).
If, on the other hand, $H\notin \mathfrak{p}$, then $I_{t-2}(\phi_h)\subset \mathfrak{p}$, hence the codimension of $\mathfrak{p}$ is $\geq n-3-(t-2)+2=n-t+1$ using (\ref{F1_of_h}).

For $t=n-1$ the result is clear.
\qed

\begin{rmk}\rm It is not apparent what would be a common natural generalization of the above examples.
Having $g=x_n^d+g_1(x_1,\ldots,x_{n-1})$ with $g_1$ being a reduced free divisor of degree $d$, while taking $h$ to be the
product of a subset of the irreducible factors of $g_1$, will not work in general.
For example, let $n=5$ and $g=x_5^5-g_1(x_1,\ldots,x_4)$, with $g_1=x_1x_2(x_1+x_2)(x_1x_4+x_2x_3)$.
Then $g_1$ is a free divisor of linear type (cf. Proposition~\ref{CalderoNarvaez2}), but adding to $g$ a product
of $3$ of the irreducible factors of $g_1$ will not result in a free divisor.

Note that the second of the above constructions can start with irreducible divisors $g,h$, with $h$ free.
Of course, iterating the above procedure by adding more variables gives other (reducible) free divisors of
quite arbitrary degrees. The simplest iterative example would initialize with $h=x_1$ and $g=x_3^2-x_2x_1$.
\end{rmk}

\subsubsection{Recap of Koszul free divisors}

The notion of a {\em Koszul free divisor} was introduced in \cite{Cal}.
Its full properties and relationship to other notions have been carried in \cite{CalNar}.
In \cite{s1} a different approach to this notion was introduced,  retaining the
original terminology.
The first approach was given in terms
of differential operators (ideal of ``differential symbols''), whereas the version in \cite{s1}
is more elementary while dealing directly with the symmetric algebra (``symmetric symbols'').
Although the two notions were given in different contexts -- one analytic, the other polynomial --
a more thorough comparison would be welcome.
In this work we will refer to the Koszul freeness according to the second approach,
hoping no confusion is expected as our setup is always polynomial and not local analytic.

According to the latter it follows that, in the (global) Eulerian situation, if the
symmetric algebra of the gradient ideal is a Cohen--Macaulay ring then saying that the divisor is Koszul free
is tantamount to freeness plus the Euler equation being
a regular element modulo the defining equations of the symmetric algebra.
This is why, in the Eulerian setup, it is fairly easy from the computational side to verify if a free divisor
is Koszul free in the present sense (see \cite[Proposition 3.11]{s1} for further details).

Let us review the latter approach in brief.

Let $L$ be a free $R$-module of finite rank and let $\{e_1,\ldots,e_p\}$ be a basis thereof.
Then ${\mathcal S}_R(L)=R[T_1,\ldots,T_p]$, a polynomial ring on variables $T_1,\ldots, T_p$.
Any element $g\in L$ is of the form $g=a_1e_1+\cdots+a_pe_p$, for some $a_i\in R$.
So one obtains $g(1):=a_1T_1+\cdots+a_pT_p\in {\mathcal S}_R(L)$.
More generally, any submodule $Z\subset L$ induces, in this way,  an ideal $Z(1)$ generated by linear forms inside ${\mathcal S}_R(L)$.
Call $Z(1)$ the {\it polynomial idealization} of the embedding $Z\subset L$.

In our setup, $Z=D_F$ and $L={\rm Der}(R)$, with $R=\mathbb K[x,y,z]$. If $\theta=\alpha\frac{\partial}{\partial x}+\beta\frac{\partial}{\partial y} +\gamma\frac{\partial}{\partial z}$ is in $D_F$, then $\theta(F)=P_F F,$ for some $P_F\in R$. Then we can write $$\theta=\underbrace{(\alpha-\frac{1}{d}xP_F)}_{A}\frac{\partial}{\partial x}+ \underbrace{(\beta-\frac{1}{d}yP_F)}_{B}\frac{\partial}{\partial y}+ \underbrace{(\gamma-\frac{1}{d}zP_F)}_{C}\frac{\partial}{\partial z}+\frac{1}{d}P_F\theta_{E},$$ where $d=\deg(F)$
and $\theta_{E}$ denotes the Euler relation.
Observe that $(A,B,C)$ is a syzygy on $(F_x,F_y,F_z)$.

With this
$$D_F(1)=\langle xT_1+yT_2+zT_3\rangle +\langle\{AT_1+BT_2+CT_3|AF_x+BF_y+CF_z=0\}\rangle$$
is the polynomial idealization of the embedding $D_F\subset {\rm Der}(R)$.
Note that ${\mathcal S}_R(J_F/(F))={\mathcal S}_R({\rm Der}(R))/D_F(1)$, thus redirecting the properties of
the polynomial idealization to those of the symmetric algebra of the $R$-module $J_F/(F)$ (an analogous result
takes place for arbitrary divisors, replacing $J_F$ by the entire gradient ideal of $F$).

According to the definition in \cite{s1}, in the present context of plane divisors, $Y=V(F)$ is a Koszul free divisor if $Y$ is a free divisor
and the polynomial idealization  $D_F(1)\subset {\mathcal S}_R({\rm Der}(R))$ has codimension $3$.

As stated in \cite[Proposition 4.1]{s1}, a free  divisor on $\mathbb P^2$ is automatically Koszul free.
Quite generally,  if $V(F)$ is free and $J_F$ of linear type then $V(F)$ is Koszul free
 (cf. \cite[Corollary 3.12]{s1}
or   \cite[Proposition 4.1 and Remark 4.2]{s1} --
but the precedence for this result is \cite[Proposition 3.2]{CalNar}).
The converse is not true, though finding such examples in the homogeneous case are not easily guessed.

Let us  change gears by revisiting a couple of known examples  (\cite[Example 1.11 and Example 6.2]{CalNar})
in the non-homogenous case that may help to clarify the theory.
For convenience we state the relevant properties in the form of a proposition.

\begin{prop}\label{CalderoNarvaez1}
Let $f = 256z^3 - 128x^2z^2 + 16x^4z + 144xy^2z - 4x^3y^2 - 27y^4\in R=\C[x,y,z]$.
Then:
\begin{enumerate}
\item[{\rm (i)}] $f$ is a weighted homogenous Koszul free divisor with  Euler equality
$12f=2xf_x+3yf_y+4zf_z$.
\item[{\rm (ii)}] The gradient ideal $I\subset R$ of $f$
is an ideal of linear type.
\item[{\rm (iii)}] Let $F\in \mathbb K[x,y,z,t]$ denote the homogenization of $f$ relative to a new variable $t$.
Then
\newline {\rm(a)} $V(F)$ is a not a free divisor of ${\mathbb P}^3$.
\newline {\rm(b)} $J_F$ is of linear type.
\end{enumerate}
\end{prop}
\demo (i) By direct inspection, $J_f$ is the ideal generated by the $2$-minors of the following
matrix
\begin{equation}\label{old_matrix}
\left(
\begin{array}{cc}
-3y  &  \quad x^2-12z\\
x^2-4z   & 3xy\\
1/2xy  &   9/4y^2-4xz
\end{array}
\right)
\end{equation}
(The matrix itself was obtained with {\em Macaulay} (\cite{Macaulay})).
By \cite[Proposition 4.1]{s1}, $V(f)$ is Koszul free.
By direct verification, it is clear that $f$ is weighted homogeneous with weights $2,3,4$ for $x,y,z$,
respectively. This proves the form of the Euler relation.

(ii) It is straightforward to see that the entries of the above matrix generate an $\langle x,y,z\rangle $-primary ideal.
Therefore, we conclude from Proposition~\ref{old_thm}.

(iii)
(a) To show that $V(F)$ is not a free divisor one could of course do a simple calculation with \cite{Macaulay},
however we choose to give a conceptual argument that may be used elsewhere.

By Lemma~\ref{lin_syzygy}, $J_F$ admits the linear syzygy $(x,-3/2y, -4z, 6t)^t$.
On the other hand, the leftmost syzygy of $J_f$ in (\ref{old_matrix}) gives by homogenization
a syzygy of $J_F$ of degree $2$ which is visibly not a multiple of the linear one.
Clearly, $J_F$  cannot be generated the $2$-minors of these two syzygies (either because it is $4$-generated
or by degree reason since $\deg(F)=5$).
Therefore, it admits at least one additional minimal syzygy independent of the first two.
Showing that any of the additional ones has standard degree
at least $2$ will do it since the sum of the degrees of any $3$ minimal syzygies would exceed the degree of the
partial derivatives of $F$, which is $4$.
Equivalently, we show that $J_F$ has no additional minimal linear syzygies.

To prove the latter claim, suppose that $\ell_xF_x+\ell_yF_y+\ell_zF_z+\ell_tF_t=0$, with the $\ell$'s linear forms in $\mathbb K[x,y,z,t]$.
Mapping $t\mapsto 1$ gives the relation
\begin{equation}\label{linear_relation_of_J_F}
\widetilde{\ell_x}f_x+\widetilde{\ell_y}f_y+\widetilde{\ell_z}f_z=-\widetilde{\ell_t}\widetilde{F_t},
\end{equation}
where $\,\tilde{}\,$ indicates  the resulting polynomial in $\mathbb K[x,y,z]$.
On the other hand, dehomogenizing the Euler relation $5F=xF_x+yF_y+zF_z+tF_t$ yields
\begin{equation}\label{Euler_of_F_dehomogenized}
5f=xf_x+yf_y+zf_z+\widetilde{F_t}.
\end{equation}
Using the weighted homogeneous Euler relation $12f=2xf_x+3yf_y+4zf_z$ and (\ref{Euler_of_F_dehomogenized})
furnishes $\widetilde{F_t}=(-5/6)xf_x+(5/4)yf_y+(10/3)zf_z$.
Comparing with (\ref{linear_relation_of_J_F}) leads to the relation
$$(\widetilde{\ell_x}-(5/6)\widetilde{\ell_t}x)f_x +(\widetilde{\ell_y}+(5/4)\widetilde{\ell_t}y)f_y
+  (\widetilde{\ell_z}+(10/3)\widetilde{\ell_t}z)f_z=0.  $$
Since $f$ is weighted homogeneous, so are its partial derivatives (of weighted degrees $10, 9, 8$, respectively).
Therefore the coefficients of the last relation must be weighted homogeneous as well.
But since a polynomial of degree $1$ can only be weighted homogeneous for distinct weights of the variables
(as is the case here) if it is a $k$-multiple of a variable, a close scrutiny of the above coefficients
shows that they are not weighted homogeneous unless they vanish.
In this case, one gets the following equalities
$$\left\{
\begin{array}{l}
\widetilde{\ell_x}=(5/6)\widetilde{\ell_t}x\\
\widetilde{\ell_y}=-(5/4)\widetilde{\ell_t}y\\
\widetilde{\ell_z}=-(10/3)\widetilde{\ell_t}z
\end{array}
\right.
$$
Since these are equalities of polynomials with standard degrees, we must have $\widetilde{\ell_t}=\alpha$,
for some $\alpha\in k$. It follows that $\widetilde{\ell_t}=\alpha t$, hence
$$(\ell_x,\ell_y,\ell_z,\ell_t)^t=\frac{6}{5}\,\alpha\,(x,-3/2y, -4z, 6t)^t,$$
showing that the assumed linear syzygy is a multiple of the linear syzygy already obtained.

\smallskip

(b) To show that $J_F$ is an ideal of linear type, one uses that it is an ideal of deviation $2$.
For deviation $2$ Cohen--Macaulay ideals  in a Cohen--Macaulay ring, it is known that their Koszul homology modules are Cohen--Macaulay
(\cite{AvHe}).
Therefore, by \cite[Theorem 9.1]{HSV1} it suffices to show that $J_F$ is locally at every prime ideal $P$ generated
by at most $\dim R_P$ elements.
This is obvious if $P$ is the maximal homogeneous ideal of $\mathbb K[x,y,z,t]$.
If $\dim R_P=3$, the linear syzygy of $J_F$ shows that $(J_F)_P$ is generated by at most $3$ elements.

It remains to prove that $J_F$ is locally a complete intersection at its minimal primes.
Now, the singular points of $V(F)$ in $\mathbb{P}^2$ are either those of $V(f)$ in $\mathbb{A}^2$ ($t=1$)
or else lie on the plane $t=0$.
A computation with \cite{Macaulay} yields that the radical of $J_f$ is generated by the polynomials
$$x^3+9/2y^2-4xz,\, y(x^2+12z),\, xy^2-8/3x^2z+32/3z^2.$$
Therefore, any minimal prime of $J_f$ contains either $y$ or $x^2+12z$.
From this one gets that the minimal primes are $P:=\langle x^2-4z,y\rangle $ and $Q:=\langle x^2+12z, y^2-32/9xz\rangle $.
Now, the entry $x^2-12z$ of the matrix (\ref{old_matrix}) does not belong to $P$, hence
locally at $P$ the first partial derivative belongs to the ideal generated by the other two partial derivatives.
Similarly, the entry $-3y$ does not belong to $Q$, hence locally at $Q$ the ideal $J_f$ is again a complete intersection.

It remains to check the minimal primes of $J_F$ that contain $t$.
Another computation with \cite{Macaulay} gives that the radical of $J_F$ is generated by the homogenized
generators of $J_f$ and the additional polynomial $y^3-32/9xyz$.
It follows from this that $t$ does not belong to any minimal prime of $J_F$ (i.e., $V(F)$ has no singularities at
infinity).

This finishes the proof.
\qed

\begin{prop}\label{CalderoNarvaez2}
Let $f = xy(x + y)(x + yz)\in R=\C[x,y,z]$ and let $J_f\subset R$ denote the ideal generated by its partial derivatives.
Then:
\begin{enumerate}
\item[{\rm (i)}] $f$ is a reduced free Eulerian divisor, but not quasihomogeneous in the sense
of all weights positive  -- namely, one has $f=1/4xf_x+1/4yf_y$.
\item[{\rm (ii)}] $J_f$ is not an ideal of linear type.
\item[{\rm (iii)}] The symmetric algebra ${\mathcal S}_R(J_f)$ is a complete intersection on
which the Euler equation $e$ is a zero-divisor  -- in particular, $f$ is not
Koszul free.
\item[{\rm (iv)}] Letting $F\in R[t]$ denote the homogenization of $f$ relative to a variable $t$, one has:
\newline {\rm (a)} $F$ is a free divisor of $\mathbb{P}^3$, but not Koszul free.
\newline {\rm (b)} The symmetric algebra of $J_F$ is a reduced Cohen--Macaulay ring,
and $J_F$ is not syzygetic.
\newline {\rm (c)} $F$ is homaloidal.
\end{enumerate}
\end{prop}
\demo (i) The Euler-like expression of $f$  is a direct calculation with the partial derivatives.
To see that $f$ is a free divisor we compute with \cite{Macaulay}, finding that the syzygy matrix
of $J_f$ is
\begin{equation}\label{old_matrix}
\left(
\begin{array}{cc}
xy  &  \quad x^2+4xy\\
y^2   & -3xy\\
-4x-4yz  &   -8x+4xz-4yz
\end{array}
\right)
\end{equation}

(ii) As one easily checks, the ideal generated by the entries of the matrix is contained in $\langle x,y\rangle$.
Therefore, $J_f$ is not of linear type.

(iii) The assertion about ${\mathcal S}_R(J_f)$ being a complete intersection is clear.
By the same token as in (ii), $\langle x,y\rangle R[T,U,V]$ is a minimal prime of the ideal $\langle \mathcal{J}_1\rangle $
of definition of ${\mathcal S}_R(J_f)$
over the polynomial ring $R[T,U,V]$.
Since the Euler relation read in degree one is $e:=1/4(xT+yU)\in \langle x,y\rangle R[T,U,V]$, it is a zerodivisor on ${\mathcal S}_R(J_f)$.

(iv) Nearly everything in this item needs computer intervention, such as with \cite{Macaulay}.
A highlight is the structural matrix of $J_F$:
$$
\left(
\begin{array}{ccc}
0  & x & 0\\
0 & y   & xy+y^2\\
x & 4z  &   -3yz\\
-y & -4t & -2yz-xt-2yt
\end{array}
\right),
$$
where the linear syzygies are easily guessed without computer calculation, due to the particular shape of $f$ and its
``nearly'' weighted homogeneous structure.
For the first linear syzygy, note $z$ and $t$ have symmetric roles.
In particular, one has $\partial F/\partial z=xy^2(x+y)$
and $\partial F/\partial t=x^2y(x+y)$.
The second syzygy is like the one obtained in Lemma~\ref{lin_syzygy}, pretending that the Euler relation is a
weighted homogenous relation.

In item (b), reducedness and minimal primes are not easily seen without resorting to the computer.
Therefore, nor is the lack of Koszul freeness.
An additional feature is that, in analogy to the case of $J_f$, saturating the defining ideal of the symmetric algebra
of $J_F$ by the (homogeneous) idealized Euler relation yields the defining ideal of the Rees algebra of $J_F$.
The latter again acquires an additional generator of bidegree $(1,2)$.
Together with the two linear syzygies, it produces enough relations to conclude by \cite[Theorem 2.18]{AHA} that the
weak Jacobian dual matrix has rank $3$, and hence the rational map defined by the partial derivatives of $J_F$
is a Cremona map of $\mathbb{P}^3$.
This gives (c).
\qed

\begin{rmk}\rm
In part (iii) of the above one can be a lot more prolixious. Namely, writing
\begin{eqnarray*}
{\mathcal S}_R(J_f)/(f)&\simeq &{\mathcal S}_R({\rm Der}(R))/D_f(1)\simeq  R[T,U,V]/{\mathcal D}\\
{\mathcal S}_R(J_f)&\simeq &R[T,U,V]/\langle {\mathcal J}_1\rangle \surjects {\mathcal R}_R(J_f)\simeq  R[T,U,V]/{\mathcal J},
\end{eqnarray*}
one has
$$\langle {\mathcal J}_1\rangle ={\mathcal J}\cap {\mathcal D}, \quad {\mathcal D}= \langle {\mathcal J}_1, e\rangle \quad {\rm and}\quad
{\mathcal J}={\mathcal J}_1:e=\langle {\mathcal J}_1, q\rangle ,$$
for a suitable $q\notin \langle {\mathcal J}_1\rangle $ of degree $2$ in $T,U,V$.
In particular, $J_f$ is not even syzygetic.
Finally, ${\mathcal R}_R(J_f)$ is Cohen--Macaulay.

The details of part (iv) are pretty much similar, except for the reducedness of the symmetric algebra --
a matter of wondering.
\end{rmk}

Clearly, $f$ is here a fairly degenerated example, sharing a proper factor with one of its derivatives -- something
that cannot happen if $f$ is irreducible.
At the other end, homaloidal (homogeneous) polynomials are also quite rare, so one wonders if there is any
bridge to freeness under suitable conditions.

In the light of the above example, it seems natural to ask:

\begin{quest}\rm
Let $f\in R=k [x_1,\ldots,x_n]$ be a reduced Eulerian divisor and let $F\in R[t]$ denote its homogenization.
If $f$ is free (respectively, Koszul free), when is $F$ free (respectively, Koszul free)?
\end{quest}

If one does not assume the Euler condition then there is a huge class of counter-examples.
Namely, take a homogeneous  irreducible $F\in  \mathbb K[x_1,\ldots,x_{n+1}] (n\geq 2)$ whose associated projective
hypersurface  is smooth and let $f\in \mathbb K[x_1,\ldots,x_n]$ denote one of its dehomogenizations.
This is because the partial derivatives of $F$ generate a
complete intersection  of codimension $n+1\geq 3$.
On the other hand, if the hypersurface defined by $f$ is smooth but its projective closure has singular points
then the issue remains (see Proposition~\ref{free_quintics}  where dehomogenization $f$ at $y=1$
is smooth but the projective closure has a singular point elsewhere -- here we have an affirmative answer to
the above question for such $f$).

\renewcommand{\baselinestretch}{1.0}
\small\normalsize 

\bibliographystyle{amsalpha}

\end{document}